\documentclass[a4paper,11pt]{amsart}

\usepackage{amsmath,amsxtra,amssymb,latexsym, amscd,amsthm,dsfont}
\usepackage[mathscr]{eucal}
\usepackage{mathrsfs,bbm,enumerate,calc,yhmath,mathtools,faktor}
\usepackage{stackengine}
\usepackage{pict2e}
\usepackage{tikz}
\usepackage{graphicx}
\usepackage{color}
\usepackage[foot]{amsaddr}

\usepackage[a4paper]{geometry}
\usepackage{hyperref}
\hypersetup{
colorlinks=true, 
breaklinks=true, 
urlcolor= blue, 
linkcolor= black, 
citecolor=black, 
pdftitle={Exposition}, 
pdfauthor={Nalini Anantharaman and Laura Monk}, 
pdfsubject={} 
}
\usepackage[capitalize]{cleveref}
\crefname{equation}{equation}{equations}
\Crefname{equation}{Equation}{Equations}

\numberwithin{equation}{section}

\title[Spectral gap of random hyperbolic surfaces]{Spectral gap of random hyperbolic surfaces}

\author{Nalini Anantharaman\textsuperscript{1} and Laura Monk\textsuperscript{2}}

\address[1]{Coll\`ege de France, 11 place Marcelin Berthelot, 75005 Paris,
  France} 

\address[2]{School of Mathematics, University of Bristol, Bristol BS8 1UG, U.K.}

\email{nalini.anantharaman@college-de-france.fr}
\email{laura.monk@bristol.ac.uk}
 
\makeatletter
\@namedef{subjclassname@2020}{\textup{2020} Mathematics Subject Classification}
\makeatother

\subjclass[2020]{Primary 58J50, 32G15; Secondary 05C80, 11F72}

\keywords{Random hyperbolic surfaces, Weil--Petersson form, moduli space,
  spectral gap, closed geodesic, Selberg trace formula.}

\theoremstyle{plain}
\newtheorem{thm}{Theorem}[section]
\newtheorem{prp}[thm]{Proposition}
\newtheorem{cor}[thm]{Corollary}
\newtheorem{lem}[thm]{Lemma}

\newtheorem*{namedthm}{\namedthmname}
\newcounter{namedthm}



\theoremstyle{definition}
\newtheorem{defa}[thm]{Definition}

\newtheorem{rem}[thm]{Remark}

\newtheorem{nota}[thm]{Notation}
\newtheorem*{summ}{Summary}

\renewcommand{\d}{\, \mathrm{d}}

\makeatletter
\newcommand*{\ov}[1]{%
  $\m@th\overline{\mbox{#1}}$%
}
\newcommand*{\ovA}[1]{%
  $\m@th\overline{\mbox{#1}\raisebox{3mm}{}}$%
}
\newcommand*{\ovB}[1]{%
  $\m@th\overline{\mbox{#1\rule{0pt}{3mm}}}$%
}
\newcommand*{\ovC}[1]{%
  $\m@th\overline{\mbox{#1\strut}}$%
}
\newcommand*{\ovD}[1]{%
  $\m@th\overline{\mbox{#1\vphantom{\"A}}}$%
}
\newcommand*{\ovE}[1]{%
  $\m@th\overline{\raisebox{0pt}[1.2\height]{#1}}$%
}
\newcommand*{\ovF}[1]{%
  $\m@th\overline{\raisebox{0pt}[\dimexpr\height+0.3mm\relax]{#1}}$%
}
\newcommand*{\ovG}[1]{%
  $\m@th\overline{\raisebox{0pt}[\dimexpr\height+1mm\relax]{#1\vphantom{A}}}$%
}
\makeatother

\newcommand\runderset[2][\sim]{\mathrel{\ensurestackMath{%
  \stackengine{-.2pt}{\scriptscriptstyle#2}{#1}{O}{c}{F}{F}{S}}}}
\newcommand{\N}{\mathbb{N}}

\newcommand{\R}{\mathbb{R}}
\newcommand{\C}{\mathbb{C}}

\newcommand{\D}{\mathcal{D}}

\DeclareMathOperator{\Id}{id}

\DeclareMathOperator{\hyp}{hyp}

\DeclareMathOperator\argcosh{argch}
\DeclareMathOperator\argsh{argsh}

\DeclareSymbolFont{extraup}{U}{zavm}{m}{n}
\DeclareMathSymbol{\varheart}{\mathalpha}{extraup}{86}
\DeclareMathSymbol{\vardiamond}{\mathalpha}{extraup}{87}

\newcommand{\nwc}{\newcommand}

\nwc{\mf}{\mathbf} 
\nwc{\blds}{\boldsymbol} 
\nwc{\ml}{\mathcal} 

\nwc{\lam}{\lambda}
\nwc{\del}{\delta}
\nwc{\Del}{\Delta}
\nwc{\Lam}{\Lambda}
\nwc{\elll}{\ell}

\nwc{\IA}{\mathbb{A}} 
\nwc{\IB}{\mathbb{B}} 
\nwc{\IC}{\mathbb{C}} 
\nwc{\ID}{\mathbb{D}} 
\nwc{\IE}{\mathbb{E}} 
\nwc{\IF}{\mathbb{F}} 
\nwc{\IG}{\mathbb{G}} 
\nwc{\IH}{\mathbb{H}} 
\nwc{\IN}{\mathbb{N}} 
\nwc{\IP}{\mathbb{P}} 
\nwc{\IQ}{\mathbb{Q}} 
\nwc{\IR}{\mathbb{R}} 
\nwc{\IS}{\mathbb{S}} 
\nwc{\IT}{\mathbb{T}} 
\nwc{\IZ}{\mathbb{Z}} 
\def\bbbone{{\mathchoice {1\mskip-4mu {\rm{l}}} {1\mskip-4mu {\rm{l}}}
{ 1\mskip-4.5mu {\rm{l}}} { 1\mskip-5mu {\rm{l}}}}}
\def\bbleft{{\mathchoice {[\mskip-3mu {[}} {[\mskip-3mu {[}}{[\mskip-4mu {[}}{[\mskip-5mu {[}}}}
\def\bbright{{\mathchoice {]\mskip-3mu {]}} {]\mskip-3mu {]}}{]\mskip-4mu {]}}{]\mskip-5mu {]}}}}
\nwc{\setK}{\bbleft 1,K \bbright}
\nwc{\setN}{\bbleft 1,\cN \bbright}
 \newcommand{\Lim}{\mathop{\longrightarrow}\limits}

\nwc{\va}{{\bf a}}
\nwc{\vb}{{\bf b}}
\nwc{\vc}{{\bf c}}
\nwc{\vd}{{\bf d}}
\nwc{\ve}{{\bf e}}
\nwc{\vf}{{\bf f}}
\nwc{\vg}{{\bf g}}
\nwc{\vh}{{\bf h}}
\nwc{\vi}{{\bf i}}
\nwc{\vI}{{\bf I}}
\nwc{\vj}{{\bf j}}
\nwc{\vk}{{\bf k}}
\nwc{\vl}{{\bf l}}
\nwc{\vm}{{\bf m}}
\nwc{\vM}{{\bf M}}
\nwc{\vn}{{\bf n}}
\nwc{\vo}{{\it o}}
\nwc{\vp}{{\bf p}}
\nwc{\vq}{{\bf q}}
\nwc{\vr}{{\bf r}}
\nwc{\vs}{{\bf s}}
\nwc{\vt}{{\bf t}}
\nwc{\vu}{{\bf u}}
\nwc{\vv}{{\bf v}}
\nwc{\vw}{{\bf w}}
\nwc{\vx}{{\bf x}}
\nwc{\vy}{{\bf y}}
\nwc{\vz}{{\bf z}}
\nwc{\bal}{\blds{\alpha}}
\nwc{\bep}{\blds{\epsilon}}
\nwc{\barbep}{\overline{\blds{\epsilon}}}
\nwc{\bnu}{\blds{\nu}}
\nwc{\bmu}{\blds{\mu}}
\nwc{\bet}{\blds{\eta}}
\nwc{\rK}{\mathrm{K}}
\nwc{\n}{n_{\mathbf{S}}}
\nwc{\g}{g_{\mathbf{S}}}

\nwc{\bk}{\blds{k}}
\nwc{\bm}{\blds{m}}
\nwc{\bM}{\blds{M}}
\nwc{\bp}{\blds{p}}
\nwc{\bq}{\blds{q}}
\nwc{\bn}{\blds{n}}
\nwc{\bv}{\blds{v}}
\nwc{\bw}{\blds{w}}
\nwc{\bx}{\blds{x}}
\nwc{\bxi}{\blds{\xi}}
\nwc{\by}{\blds{y}}
\nwc{\bz}{\blds{z}}

\newcommand{\tfL}{\omega}

\nwc{\cA}{\ml{A}}
\nwc{\cB}{\ml{B}}
\nwc{\cC}{\ml{C}}
\nwc{\cD}{\ml{D}}
\nwc{\cE}{\ml{E}}
\nwc{\cF}{\ml{F}}
\nwc{\cG}{\ml{G}}
\nwc{\cH}{\ml{H}}
\nwc{\cI}{\ml{I}}
\nwc{\cJ}{\ml{J}}
\nwc{\cK}{\ml{K}}
\nwc{\cL}{\ml{L}}
\nwc{\cM}{\ml{M}}
\nwc{\cN}{\ml{N}}
\nwc{\cO}{\ml{O}}
\nwc{\cP}{\ml{P}}
\nwc{\cQ}{\ml{Q}}
\nwc{\cR}{\ml{R}}
\nwc{\cS}{\ml{S}}
\nwc{\cT}{\ml{T}}
\nwc{\cU}{\ml{U}}
\nwc{\cV}{\ml{V}}
\nwc{\cW}{\ml{W}}
\nwc{\cX}{\ml{X}}
\nwc{\cY}{\ml{Y}}
\nwc{\cZ}{\ml{Z}}

\nwc{\fA}{\mathfrak{a}}
\nwc{\fB}{\mathfrak{b}}
\nwc{\fC}{\mathfrak{c}}
\nwc{\fD}{\mathfrak{d}}
\nwc{\fE}{\mathfrak{e}}
\nwc{\fF}{\mathfrak{f}}
\nwc{\fG}{\mathfrak{g}}
\nwc{\fH}{\mathfrak{h}}
\nwc{\fI}{\mathfrak{i}}
\nwc{\fJ}{\mathfrak{j}}
\nwc{\fK}{\mathfrak{k}}
\nwc{\fL}{\mathfrak{l}}
\nwc{\fM}{\mathfrak{m}}
\nwc{\fN}{\mathfrak{n}}
\nwc{\fO}{\mathfrak{o}}
\nwc{\fP}{\mathfrak{p}}
\nwc{\fQ}{\mathfrak{q}}
\nwc{\fR}{\mathfrak{r}}
\nwc{\fS}{\mathfrak{s}}
\nwc{\fT}{\mathfrak{t}}
\nwc{\fU}{\mathfrak{u}}
\nwc{\fV}{\mathfrak{v}}
\nwc{\fW}{\mathfrak{w}}
\nwc{\fX}{\mathfrak{x}}
\nwc{\fY}{\mathfrak{y}}
\nwc{\fZ}{\mathfrak{z}}

\nwc{\bT}{\mathbf{T}}

\nwc{\tA}{\widetilde{A}}
\nwc{\tB}{\widetilde{B}}
\nwc{\tE}{E^{\vareps}}
\nwc{\tk}{\tilde k}
\nwc{\tN}{\tilde N}
\nwc{\tP}{\widetilde{P}}
\nwc{\tQ}{\widetilde{Q}}
\nwc{\tR}{\widetilde{R}}
\nwc{\tV}{\widetilde{V}}
\nwc{\tW}{\widetilde{W}}
\nwc{\ty}{\tilde y}
\nwc{\teta}{\tilde \eta}
\nwc{\tdelta}{\tilde \delta}
\nwc{\tlambda}{\tilde \lambda}
\nwc{\ttheta}{\tilde \theta}
\nwc{\tvartheta}{\tilde \vartheta}
\nwc{\tPhi}{\widetilde \Phi}
\nwc{\tpsi}{\tilde \psi}
\nwc{\tmu}{\tilde \mu}

\nwc{\To}{\longrightarrow} 

\nwc{\ad}{\rm ad}
\nwc{\eps}{\epsilon}
\nwc{\ep}{\epsilon}
\nwc{\vareps}{\varepsilon}

\def\ep{\epsilon}

\def\sq2{\sqrt{2}}

\def\t2{{\mathbb T}^2}
\def\s2{{\mathbb S}^2}

\def\N{\mathbb{N}}

\def\R{\mathbb{R}}

\def\C{\mathbb{C}}
\def\O{\mathcal{O}}

\nwc{\lap}{\bigtriangleup}
\nwc{\rest}{\restriction}
\nwc{\Diff}{\operatorname{Diff}}
\nwc{\diam}{\operatorname{diam}}
\nwc{\Res}{\operatorname{Res}}
\nwc{\Spec}{\operatorname{Spec}}
\nwc{\Vol}{\operatorname{Vol}}
\nwc{\Op}{\operatorname{Op}}
\nwc{\supp}{\operatorname{supp}}
\nwc{\Span}{\operatorname{span}}

\nwc{\dia}{\varepsilon}
\nwc{\cut}{f}
\nwc{\qm}{u_\hbar}

\def\hto0{\xrightarrow{\hbar\to 0}}

\def\rto0{\xrightarrow{r\to 0}}

\providecommand{\norm}[1]{\lVert#1\rVert}

\nwc{\la}{\langle}
\nwc{\ra}{\rangle}
\nwc{\lp}{\left(}
\nwc{\rp}{\right)}

\nwc{\bequ}{\begin{equation}}
\nwc{\be}{\begin{equation}}
\nwc{\ben}{\begin{equation*}}
\nwc{\bea}{\begin{eqnarray}}
\nwc{\bean}{\begin{eqnarray*}}
\nwc{\bit}{\begin{itemize}}
\nwc{\bver}{\begin{verbatim}}

\nwc{\eequ}{\end{equation}}
\nwc{\ee}{\end{equation}}
\nwc{\een}{\end{equation*}}
\nwc{\eea}{\end{eqnarray}}
\nwc{\eean}{\end{eqnarray*}}
\nwc{\eit}{\end{itemize}}
\nwc{\ever}{\end{verbatim}}

\makeatletter
\newlength{\temp@wc@width}
\newlength{\temp@wc@height}
\newcommand{\widecheck}[1]{%
  \setlength{\temp@wc@width}{\widthof{$#1$}}%
  \setlength{\temp@wc@height}{\heightof{$#1$}}%
  #1\hspace{-\temp@wc@width}%
  \raisebox{\temp@wc@height+2pt}[\heightof{$\widehat{#1}$}]%
     {\rotatebox[origin=c]{180}{\vbox to 0pt{\hbox{$\widehat{\hphantom{#1}}$}}}}%
}
\makeatother

\newcommand{\Pwpo}{\mathbb{P}_g^{\mathrm{\scriptsize{WP}}}}
\newcommand{\Ewpo}[1][g]{\mathbb{E}_{#1}^{\mathrm{\scriptsize{WP}}}}

\newcommand{\Pwp}[1]{\Pwpo \left( #1 \right)}
\newcommand{\Ewp}[2][g]{\mathbb{E}_{#1}^{\mathrm{\scriptsize{WP}}} \Bigg[ #2 \Bigg]}

\DeclarePairedDelimiter{\paren}{(}{)}

\DeclarePairedDelimiter{\abso}{|}{|}
\DeclarePairedDelimiter{\brac}{[}{]}

\let\div\relax
\newcommand{\div}[1]{\paren*{\frac{#1}{2}}}

\renewcommand{\O}[2][ ]{\mathcal{O}_{#1} \left( #2 \right)}
\newcommand{\Ow}[2][ ]{\mathcal{O}_{#1}^w \left( #2 \right)}

\newcommand{\1}[1]{\mathds{1}_{#1}}

\newcommand{\av}[2][\mathrm{all}]{\langle #2 \rangle_g^{{#1}}}
\newcommand{\avN}[2][\mathrm{all}]{\langle #2 \rangle_{g, Q}^{{#1}}}
\newcommand{\BigavN}[2][\mathrm{all}]{\Big\langle #2 \Big\rangle_{g, Q}^{{#1}}}
\newcommand{\avTbtf}[2][\mathrm{all}]{\langle #2 \rangle_{g, \kappa, \tfL}^{{\type}}} 
 
\newcommand{\avb}[2][\mathrm{all}]{\left\langle #2 \right\rangle_g^{{#1}}}
\newcommand{\avbtf}[2][]{\left\langle #2 \right\rangle_{g, \kappa, \tfL, \chic}^{{#1}}}

\newcommand{\FR}{\mathcal{F}}
\newcommand{\FRrem}{\mathcal{R}}

\newcommand{\eqc}[1]{[ #1 ]_{\mathrm{loc}}}
\newcommand{\eq}{\, \raisebox{-1mm}{$\runderset{\mathrm{loc}}$} \,}

\newcommand{\MC}{\mathrm{MC}}

\newcommand{\x}{\mathbf{x}}

\makeatletter
\newcommand{\smallbullet}{} 
\DeclareRobustCommand\smallbullet{%
  \mathord{\mathpalette\smallbullet@{0.5}}%
}
\newcommand{\smallbullet@}[2]{%
  \, \vcenter{\hbox{\scalebox{#2}{$\m@th#1\bullet$}}} \,%
}
\makeatother

\newcommand{\Amin}{a}
\newcommand{\domain}{\mathfrak{D}}
\newcommand{\LocTFlc}{\mathrm{Loc}_{\Sf,\chic, \chitau}^{ \kappa, \tfL, L}}
\newcommand{\tf}{\mathrm{TF}_g}
\newcommand{\geod}{\mathcal{G}}
\newcommand{\primitive}{\mathcal{P}}
\newcommand{\Atf}{\mathrm{TF}_g^{\kappa,\tfL}}
\newcommand{\Sf}{\mathbf{S}}
\newcommand{\type}{\mathbf{T}}
\newcommand{\curve}{\mathbf{c}}
\newcommand{\chic}{{\chi_+}}
\newcommand{\chitau}{{\chi_+'}}
\newcommand{\chictau}{{\chi_+''}}

\begin{document}

\begin{abstract} 
Let $X$ be a closed, connected, oriented surface of genus $g$, with a hyperbolic metric chosen at random according to the Weil--Petersson measure on the moduli space of Riemannian metrics. Let $\lambda_1=\lambda_1(X)$ be
the first non-zero eigenvalue of the Laplacian on $X$ or, in other words, the spectral gap.
In this paper we give a full road-map to prove that for arbitrarily small~$\alpha>0$,
\begin{align*}  \Pwp{\lambda_1 \leq \frac{1}{4} - \alpha^2 } \Lim_{g\To +\infty} 0.
\end{align*}
The full proofs are deferred to separate papers.
 \end{abstract}

 \maketitle

 \tableofcontents

\section{Introduction}

Let $X$ be a closed, connected, oriented surface of genus $g$, with a hyperbolic metric chosen at
random according to the Weil--Petersson measure on the moduli space $\cM_g$ of Riemannian
metrics. This measure is known to be finite, of total mass $V_g$; we normalize it to be a
probability measure, denoted by $\Pwpo$. The aim of this article is to establish asymptotic results
true \emph{with high probability}, i.e. with probability going to $1$ in the large genus limit
$g\To +\infty$. In particular, let $\lambda_1=\lambda_1(X)$ be the first non-zero eigenvalue of the
Laplacian on $X$, known as its \emph{spectral gap}. Motivations to study this quantity can be
found in \cite{hide2021,Ours1}. In \cite{Ours1}, we proved that for $\alpha = \frac16$ and $\eps>0$
arbitrary, we have
\begin{align} \label{e:dream} \Pwp{\lambda_1 \leq \frac{1}{4} - \alpha^2 - \epsilon} \Lim_{g\To +\infty} 0.
\end{align}
In other words, we proved that, for any $\epsilon >0$, $\lambda_1 \geq \frac{2}{9} - \epsilon$ with
high probability.  Two previous independent papers due to Wu--Xue~\cite{wu2022} and
Lipnowski--Wright~\cite{lipnowski2021} proved \eqref{e:dream} for $\alpha=\frac14$, i.e. that
$\lambda_1 \geq \frac{3}{16} - \epsilon$ with high probability.

In this paper, we describe all the steps of our strategy to prove \eqref{e:dream} for arbitrarily
small $\alpha$. This allows to conclude that, for any $\epsilon >0$,
$\lambda_1 \geq \frac 14 - \epsilon$ with high probability. In light of Huber's work
\cite{huber1974} proving that $\limsup_g \sup_{X \in \cM_g} \lambda_1(X) \leq \frac 14$, this means
that typical hyperbolic surfaces of large genus have an almost optimal spectral gap.

The aim of this expository article is to provide a detailed walkthrough of the proof, with precise
intermediate statements. A few explicit computations are provided for illustration purposes, but the
full proofs are deferred to other articles. As much as possible, we aim to make this presentation
standalone, with careful references to our four articles
\cite{anantharaman2022,Ours1,Ours2,Moebius-paper} which come together to prove~\eqref{e:dream}.

The key steps of the proof are the following.
\begin{itemize}
\item First, we define a notion of \emph{volume functions} associated to arbitrary topologies of
  closed geodesics, extending the \emph{volume polynomials} investigated by Mirzakhani for simple
  geodesics. We provide an expression and an asymptotic expansions in powers of $g^{-1}$ for these
  volume functions, presented in \S \ref{s:topo}.
\item Then, we prove that the coefficients appearing in these expansions satisfy the
  \emph{Friedman--Ramanujan} property, a newly defined notion related to on-average cancellations in
  the trace method. This statement, Theorem \ref{t:leviathan}, and its consequence,
  Proposition~\ref{cor:FR_implies_small}, yield the crucial argument to prove \eqref{e:dream}. The
  proof of Theorem \ref{t:leviathan} is the focus of our upcoming article \cite{Ours2}, and of the
  third section of this article.
\item Finally, we explain the necessity to discard a set of ``tangled surfaces'' of small but
  non-zero probability. We demonstrate how the exponential proliferation of topologies of closed
  geodesics in tangled surfaces is responsible for the failure of the naive trace method. We solve
  this issue by conditioning our argument on the set of \emph{tangle-free surfaces}. This
  delicate step is made possible by a new kind of Moebius inversion formula~\cite{Moebius-paper},
  and performed in \S \ref{s:tangles_moebius}.
\end{itemize}

 \subsubsection*{Acknowledgements} 
 This research has received funding from the EPSRC grant EP/W007010/1 and from the European Research
 Council (ERC) under the European Union’s Horizon 2020 research and innovation programme (Grant
 agreement No. 101096550).

\section{Walkthrough of the proof}
\label{s:walkthrough}

 \subsection{The trace method} \label{s:STF}

 The starting point of our analysis is, without surprise, the trace formula proven by Selberg in
 \cite{selberg1956}. It relates the spectrum of the Laplacian
$$\lambda_0(X)=0 < \lambda_1(X)\leq \lambda_2(X)\leq \ldots \To +\infty$$
on a closed connected oriented hyperbolic surface $X$ of genus $g$, to the lengths of all its
periodic geodesics. It reads, for a smooth even function $H : \R \rightarrow \R$,
\begin{equation}
  \label{eq:selberg}
  \sum_{j=0}^{+ \infty} \hat{H}(r_j(X))
  =  (g-1) \int_\R \hat{H}(r) \tanh (\pi r) r \d r
  + \sum_{\gamma \in \geod(X)} \sum_{k=1}^{+ \infty}
  \frac{\ell(\gamma) \, H(k \ell(\gamma))}{2 \sinh \div{k \ell(\gamma)}}
\end{equation}
where:
\begin{itemize}
\item for all $j$, $r_j(X) \in \R \cup i [- 1/2, 1/2]$ is a solution of
  $\lambda_j(X) = \frac 14 + r_j(X)^2$ -- the left-hand side of the formula is thus called the
  {\emph{spectral side}};
\item the Fourier transform $\hat{H}$ is
defined by $\hat{H}(r) = \int_{\R} H(\ell) \, e^{-i r \ell} \d \ell$; 
\item the first term on the right-hand side is the so-called \emph{topological
    term}, referring to the fact
  that this term only depends on the genus $g$ (in particular, when studying random hyperbolic
  surfaces of genus $g$, this term is deterministic);
\item $\geod(X)$ is the set of primitive oriented periodic geodesic on $X$ and $\ell(\gamma)$ stands
  for the length of the smooth curve $\gamma$ in $X$.
\end{itemize}

The sum over periodic geodesics is called the \emph{geometric term} of the trace formula. In the
formula, the integer $k$ represents the number of times the primitive geodesic is run over: the sum
$\sum_{k=2}^{+\infty}$ thus describes non-primitive periodic geodesics.

We draw the attention to the fact that non-simple geodesics appear in the geometric term. Dealing
with non-simple closed geodesics in the Selberg trace formula is one of the core challenges faced
when taking the average of the Selberg trace formula for random Weil--Petersson surfaces.

The Selberg trace formula holds for a class of ``nice'' functions $H$. For our purposes, we will
only consider functions $H$ of compact support, in which case both sums are absolutely convergent
\cite[Theorem 5.8]{bergeron2016}.  More precisely, let $H : \R \rightarrow \R_{\geq 0}$ be a fixed
smooth even function, with compact support $[-1,1]$, such that $\hat{H}$ is non-negative on
$\R \cup i[- \frac 12, \frac 12]$. For any $L \geq 1$, we shall apply the trace formula to
$H_L(\ell) := H(\frac{\ell}{L})$. The function $H_L$ has support $[-L, L]$, so that the geometric
term only involves periodic geodesics of length $\ell(\gamma)\leq L$. Note that $H_L$ still has
non-negative Fourier transform. The parameter $L$ will be taken to grow to $+\infty$ as
$g\To \infty$.

Eigenvalues $\lambda_j \leq\frac14$ correspond to purely imaginary $r_j $, which are are responsible
for terms such that $\hat{H}_L(r_j)$ grows roughly like $e^{|r_j| L}$.  For the first non-trivial
eigenvalue, corresponding to $j=1$, this is expressed in the following lemma:

\begin{lem}[{\cite[Lemma 3.10]{Ours1}}]
  \label{lem:growth_h_r1}
  Let $\alpha \in (0, \frac 12)$. For any $0 < \epsilon < \frac{1}{4} - \alpha^2$, there exists a
  constant $C_{\alpha,\epsilon} > 0$ such that, for any hyperbolic surface $X$, any $L \geq 1$,
  \begin{equation*}
    \lambda_1(X) \leq \frac 14 - \alpha^2 - \epsilon
    \quad \Longrightarrow \quad
    \hat{H}_{L}(r_1(X)) \geq C_{\alpha,\epsilon} \, e^{(\alpha + \epsilon) L}.
  \end{equation*}
\end{lem}
\begin{rem}Note the different roles played by $\alpha$ and $\eps$ in the discussion: $\alpha$ is
  some fixed real number, whereas $\eps$ is aimed to be an arbitrarily small number. As already
  mentioned, the papers~\cite{wu2022,lipnowski2021} dealt with the case $\alpha=\frac14$, and we
  considered $\alpha=\frac16$ in \cite{Ours1}. In order to prove~\eqref{e:dream}, we need to be able
  to take $\alpha$ fixed but arbitrarily small.
\end{rem}
\begin{rem}
  Note that, throughout this article, the test function $H$ is fixed once and for all, and every
  constant appearing can depend on the choice of $H$. For instance, the constant in Lemma
  \ref{lem:growth_h_r1} depends on $H$. 
\end{rem}

The Selberg trace formula holds for any hyperbolic surface $X$, and if $X$ is random, then both
sides also become random.  Lemma \ref{lem:growth_h_r1} provides us with the beginning of a strategy
to prove probabilistic lower bounds on $\lambda_1$. Suppose we want to prove that
$\lambda_1 \geq \frac{1}{4} - \alpha^2 - \epsilon$ with high probability.  First, we use
\cref{lem:growth_h_r1} to write
\begin{equation*}
  \Pwp{\lambda_1 \leq \frac{1}{4} - \alpha^2 - \epsilon}
  \leq \Pwp{\hat{H}_L(r_1) \geq C_{\alpha,\epsilon}  e^{(\alpha + \epsilon) L}}.
\end{equation*}
Then, the Markov inequality yields
\begin{equation} \label{e:Markov1}
  \Pwp{\lambda_1 \leq \frac{1}{4} - \alpha^2 - \epsilon}
  \leq \frac{\Ewpo \brac*{\hat{H}_L(r_1)}}{C_{\alpha,\epsilon} \, e^{(\alpha + \epsilon) L}}.
\end{equation}
In order to imply \eqref{e:dream}, it is thus sufficient to prove that, for some choice of
$L=L(g)\To +\infty$,
\begin{align} \label{e:mark}
\Ewpo \brac*{\hat{H}_L(r_1)} = o( e^{(\alpha + \epsilon) L}),
\end{align}
where by $o(\,)$ we mean that the ratio between the left-hand side and right-hand side converges to
$0$ as $g \To \infty$.  It is in our interest to take $L$ as small as possible, because we are
summing over periodic geodesics of lengths in $[0, L]$, and the number of such geodesics is known to
grow exponentially in $L$. At this stage of the discussion, $L(g)$ can grow arbitrarily slowly, but
will be forced to grow at a certain rate in the coming lines.

To control the left-hand side of \eqref{e:mark}, we can use the Selberg trace formula
\eqref{eq:selberg} and the positivity of $\hat{H}_L$, to write
\begin{align*}
\Ewpo \brac*{\hat{H}_L(r_1)}\leq
 (g-1) \int_\R \hat{H}_L(r) \tanh (\pi r) r \d r
  + \Ewpo \brac*{\sum_{\gamma \in \geod(X)} \sum_{k=1}^{+ \infty}
  \frac{\ell(\gamma) \, H_L(k \ell(\gamma))}{2 \sinh \div{k \ell(\gamma)}}}.
\end{align*}

Let us acknowledge the necessary presence of the deterministic topological term, growing linearly in
$g$ on the right-hand side: if we are to prove \eqref{e:mark}, this term forces us to take
$L\geq \frac{\log g}{\alpha+\epsilon}$.  Hence our discussion always takes place at logarithmic
scales in $g$. The crux of the analysis lies in the choice of the multiplicative constant: the
smaller $\alpha$ we aim at, the larger multiplicative constant we need.  In particular, $\alpha=0$
requires to take $L\geq \frac{\log g}{\epsilon}$ for any arbitrary $\epsilon >0$.  With this in
mind, we take from now on $L=A\log g$ with $A$ fixed.

Again because of the unavoidable presence of a term growing linearly in $g$, any contribution
bounded by a constant multiple of $g (\log g)^c$ can be grouped with the topological term and
considered to be an error term. For instance, the forthcoming Lemma \ref{lem:bound_r1} shows that
the sum $\sum_{k=2}^{+ \infty}$ in the geometric term (corresponding to non-primitive periodic
geodesics) gives such a contribution and, similarly, that the cost of replacing the factor
$\left(2\sinh (\ell(\gamma)/2)\right)^{-1}$ by $\exp(- \ell(\gamma)/2)$ is linear in $g$.

\begin{nota}
  Let us define, for any bounded compactly supported function $F$,
  \begin{equation}
    \label{e:def_av_all}
    \avb{F} := \Ewp{\sum_{\gamma \in \geod(X)} F(\ell(\gamma))}.
  \end{equation}
\end{nota}

\begin{lem}[{\cite[Lemma 3.11]{Ours1}}] 
  \label{lem:bound_r1}
   Let $L \geq 1$. If $F : \R \rightarrow \R$ is a smooth even function, supported
  on $[-L,L]$, with $\hat{F} \geq 0$ on $\R \cup i [- \frac 12, \frac
  12]$, then, for any $g\geq 2$,
  \begin{align*}
    \Ewpo \brac*{\hat{F}(r_1(X))}
  &  \leq 
      \avb{ \frac{\ell \,F(\ell)}{2 \sinh \div{\ell}} }
      + C_{F} \, L^2 g,
       \end{align*}
       for a constant $C_{F} := c \norm{F}_\infty + \norm{r \hat{F}(r)}_{\infty}< + \infty$, where
       $c$ is a universal constant independent of $g$, $L$ and $F$. The same holds replacing the
       average on the r.h.s. by $ \avb{\ell \, F(\ell) \, e^{- \frac \ell 2}}$.
\end{lem}
The proof of the lemma relies on a standard deterministic counting estimate for periodic geodesics,
that we recall here since it will be reused later in the paper. Note that since the Euler
characteristic always appears in absolute value in our estimates, we shall refer to this value as
the \emph{absolute Euler characteristic}.
\begin{lem}[{\cite[Lemma 2.1]{Ours1}}, adapted from {\cite[Theorem 4.1.6 and Lemma
  6.6.4]{buser1992}}]
  \label{lem:bound_number_closed_geod}
  Let $X$ be a hyperbolic surface, compact or bordered. For any $L>0$,
  \begin{equation*}
    \# \{ \gamma \in \geod(X) \, | \, \ell(\gamma) \leq L \}
    \leq 205 \, {\chi(X)} \, e^{L}
  \end{equation*}
  where $\chi(X)$ is the absolute Euler characteristic of $X$.  As a consequence, if $F$ is a
  function supported in $[-L, L]$,
  \begin{align}\label{e:basic}
  \sum_{\gamma \in \geod(X)}| F(\ell(\gamma))| \leq 560 \, {\chi(X)} L \norm{F(\ell)e^{\ell}}_\infty.
  \end{align}  
\end{lem}

We can use Lemma \ref{lem:bound_r1} together with \eqref{e:Markov1} to obtain that, for $L = A \log(g)$,
\begin{equation}
  \label{eq:trace_method_before_canc}
  \Pwp{\lambda_1 \leq \frac{1}{4} - \alpha^2 - \epsilon}
  = \O[\alpha, \epsilon, A]{
    \frac{\avb{\ell H_L(\ell) \, e^{- \frac \ell 2}}}{g^{(\alpha + \epsilon) A}}
    + (\log g)^2 g^{1-(\alpha+\epsilon) A}},
\end{equation}
where by $\O[x]{\,}$ we mean that the ratio between the left-hand side and the right-hand side is
bounded by a constant depending only on the parameters in the subscript $x$.

\begin{summ}
  We have shown that one can reduce the proof of \eqref{e:dream}, for a $0 < \alpha < \frac12$ and
  $0 < \epsilon < \frac{1}{4} - \alpha^2$, to the question of showing that
  \begin{align}\label{e:hope}
    \avb{ \frac{\ell \,H_L(\ell)}{2 \sinh \div{\ell}} }
    =o( e^{(\alpha +  \epsilon) L})
    \qquad \text{or} \qquad
    \avb{\ell H_L(\ell) \, e^{- \frac \ell 2}}
    =o( e^{(\alpha +  \epsilon) L})
  \end{align}
  for the length parameter $L= A \log g$, with $A\geq \frac{1}{\alpha+\epsilon}$.
\end{summ}

\subsection{The average contribution of simple geodesics}\label{s:simple}
To examine the possibility of proving~\eqref{e:hope}, let us first analyse the contribution of
simple geodesics in the average $\av{F}$, that is, geodesics without any self-intersections. This is
made feasible by Mirzakhani's integration formula. We define
\begin{equation}
  \label{eq:av_S}
  \av[\mathrm{s}]{F} := 
  \Ewp{\sum_{\gamma \text{ simple}} F(\ell(\gamma))},
\end{equation}
a similar expression to \eqref{e:def_av_all}, but where we only sum over simple oriented geodesics.
Mirzakhani's work \cite{mirzakhani2007} provides the existence of a density
$V_g^{\mathrm{s}} : \R_{>0} \rightarrow \R$ such that, for any measurable function~$F$ with compact
support,
\begin{equation}\label{e:simpledens}
  \av[\mathrm{s}]{F} = 
    \frac{1}{V_g} \int_0^{+ \infty} F(\ell) V_g^{\mathrm{s}}(\ell) \d \ell.
\end{equation}
Most importantly, Mirzakhani showed that $\ell\mapsto V_g^{\mathrm{s}}(\ell)$ is a polynomial
function (called ``volume polynomial''), and provided some recursion formulas allowing to determine
$V_g^{\mathrm{s}}$ by induction on the integer~$g$.  It follows from work of Mirzakhani and Mirzakhani--Petri
\cite{mirzakhani2013, mirzakhani2019} that there exists a constant~$c$ such that
\begin{align} \frac{V_g^{\mathrm{s}}(\ell)}{V_g}=  
    \frac{4}{\ell} \sinh^2 \div{\ell} + \O{ \frac{(1+\ell)^{c}e^{\ell}}{g}} \label{e:asympvol}. \end{align}
As a consequence,  
    \begin{align}
\nonumber  \left\langle \frac{\ell \, H_L(\ell)}{2 \sinh \div{\ell}} \right\rangle_g^{\mathrm{s}}
  & = \int_{0}^{+ \infty} \frac{\ell \, H_L(\ell)}{2 \sinh \div{\ell}} \,
 \nonumber   \frac{4}{\ell} \sinh^2 \div{\ell} \d \ell \,
+ \O{\frac{L^{c} e^{\frac L 2}}{g}}\\
   & \label{e:scont}=  2 \int_{0}^{+ \infty} H_L(\ell) \cosh \div{\ell} \d \ell
  + \O{1+\frac{L^{c} e^{\frac L 2}}{g}},
\end{align}
which is greater than $C e^{(\frac12 -\eps) L}$ for any $\eps>0$ if $[-1, 1]=\supp H$, apparently ruining our hope to prove \eqref{e:hope}.

The reason for this issue is clearly identified: we omitted to discuss the contribution of the trivial eigenvalue $\lambda_0=0$, corresponding to $r_0=\frac{i}2$, in the spectral side of the Selberg trace formula.
The trivial eigenvalue gives rise to a term growing exponentially in $L$ in the spectral side. In fact,
 \begin{align*}
 \hat{H}_L(r_0)=\hat{H}_L(i/2)= 2 \int_{0}^{+ \infty} H_L(\ell) \cosh \div{\ell} \d \ell
 \end{align*}
precisely matching the leading contribution of simple geodesics we just calculated in \eqref{e:scont}! This term appearing on both sides of the trace formula can be cancelled out, yielding
\begin{align}  \Ewpo \brac*{\hat{H}_L(r_1(X))} = \O{ L^2 g  + \frac{L^{c} e^{\frac L 2}}{g}}+ \parbox{3.5cm}{
  \begin{center}contribution of \\ non-simple geodesics\end{center}}.
\end{align}
This was observed independently by the two teams Wu--Xue and Lipnowski--Wright in \cite{wu2022,
  lipnowski2021}. Following two complementary approaches, they then remarkably achieved to show that
the contribution of non-simple geodesics can be bounded by $\O{e^{(1+\epsilon)L/2}/g} $ --
this part being, by far, the main challenge of the proof of \eqref{e:dream} for $\alpha = \frac
14$. Balancing the terms of order $e^{L/2}/g $ with the terms of order $g$ leads to the choice of
$L=4 \log g$, corresponding to $\alpha=\frac14$, yielding
 \begin{align*}\Pwp{\lambda_1 \leq  \frac{3}{16} - \epsilon} = \Pwp{\lambda_1 \leq \frac{1}{4} -\frac{1}{4^2} - \epsilon} \To 0. \end{align*}
 
 \begin{summ}
   To go further and lower the value of $\alpha$, we need to:
   \begin{itemize}
   \item understand more precisely the error terms in the asymptotics \eqref{e:asympvol} of volume
     polynomials: we want to prove asymptotic expansions to any order in powers of $g^{-1}$, and
     understand the coefficients in these expansions;
   \item extend the previous point (that is, the existence of a density \eqref{e:simpledens} and the
     asymptotics \eqref{e:asympvol}) to non-simple geodesics;
   \item find a systematic way of cancelling exponential behaviour on both sides of the averaged
     trace formula.
   \end{itemize}
   The first two points are the subject of our paper \cite{Ours1}; the results are summarized in the
   next sections.  However, it is not reasonable to expect we will be lucky enough to see, by simple
   inspection, cancellations between the spectral and geometric terms at any order in our asymptotic
   expansions.
 \end{summ}

 \subsection{Killing  dominant exponential terms}\label{s:killing}
 We propose a new way of cancelling the exponential terms of size $e^{L/2}$, that we first explore
 in the case of simple geodesics, before coming to arbitrary topologies. This strategy has first
 been implemented in our work~\cite{Ours1} in the case $\alpha=\frac16$.
 
 In order to kill the contribution of the trivial eigenvalue $0$, let us choose in the trace formula
 a test function $H_L$ whose Fourier transform $\hat H_L(r)$ vanishes if $\frac14+r^2=0$.  We simply
 do this by replacing $\hat H_L(r)$ by $(\frac14+r^2)^m\hat H_L(r)$ for some $m\geq 1$. Note that
 $r\mapsto (\frac14+r^2)^m\hat H_L(r)$ is still a non-negative function.  Equivalently, we replace
 the test function $H_L$ by $\cD^m H_L$, where $\cD$ is the differential operator
 $\cD:= \frac14-\partial^2$.  Now $\cD^m H_L$ is no longer of constant sign; this is cause for
 trouble, making necessary our Moebius inversion formula (see Remark \ref{rem:positivity}).

 Inequality \eqref{eq:trace_method_before_canc} can be replaced by:
\begin{lem}[{\cite[Lemma 3.14]{Ours1}}]
  \label{lem:selberg_reformulated}
  For any $0 < \alpha < \frac 12$, $0 < \epsilon < \frac 1 4 - \alpha^2$, $\delta > 0$, $A \geq 1$
  and $m \geq 1$, there exists a constant $C = C(H,\alpha, \epsilon, \delta, A, m)$ such that, for
  any large enough integer $g$ and for the length-scale $L = A \log g$,
  \begin{align}
    \Pwp{\delta \leq \lambda_1 \leq \frac 1 4 - \alpha^2 - \epsilon}
   \nonumber  &\leq 
     \frac{C}{g^{(\alpha + \epsilon) A}} \left(
    \avb{\frac{\ell \, \cD^m H_L(\ell)}{2 \sinh \div{\ell}}}
    + \,{g \log(g)^2}\right),
     \end{align}
and the same holds replacing the average on the r.h.s. by $\avb{\ell \,e^{-\ell/2} \cD^m H_L(\ell)}$.
\end{lem}

\begin{rem} The additional constraint $\delta \leq \lambda_1$ is necessary to deal with the fact that the modified test function $(\frac14+r^2)^m\hat H_L(r)$ vanishes for $\frac14+r^2=0$, hence it cannot be used to
exclude the possibility of arbitrarily small eigenvalues.
However, thanks to Mirzakhani's work~\cite{mirzakhani2013}, we know that there exists $\delta>0$ such that
 \begin{align*}
    \Pwp{ \lambda_1 >\delta}  \Lim_{g\To +\infty} 0.
    \end{align*}
    Thus, in order to prove \eqref{e:dream}, it is sufficient to prove that 
     \begin{align*}
    \Pwp{\delta \leq \lambda_1 \leq \frac 1 4 - \alpha^2 - \epsilon}  \Lim_{g\To +\infty} 0.
    \end{align*}
\end{rem}

Take $m=1$. Looking back at the calculation \eqref{e:scont} with $H_L$ replaced by $\cD H_L$, we now find
     \begin{align}
\nonumber  \left\langle \frac{\ell \, \cD H_L(\ell)}{2 \sinh \div{\ell}} \right\rangle_g^{\mathrm{s}}
  & = \int_{0}^{+ \infty} \frac{\ell \, \cD H_L(\ell)}{2 \sinh \div{\ell}} \,
 \nonumber   \frac{4}{\ell} \sinh^2 \div{\ell} \d \ell
  + \O{\frac{L^{c}  e^{\frac L 2}}{g}} \\  
   & \nonumber= 2 \int_{0}^{+ \infty}  \cD H_L(\ell)  \sinh \div{\ell}  \d \ell
  + \O{\frac{L^{c} e^{\frac L 2}}{g}} \\
  &\label{e:ibp} \overset{\text{IBP}}{=}
    2 \int_{0}^{+ \infty}  H_L(\ell)  \, \cD \left[\sinh \div{\ell} \right] \d \ell - H_L(0)
  + \O{\frac{L^{c}  e^{\frac L 2}}{g}} \\
  &\label{e:vanish}= 0 - H_L(0)
  + \O{\frac{L^{c}  e^{\frac L 2}}{g}}
 \end{align}
 where \eqref{e:ibp} relies on integration by parts. In equation \eqref{e:vanish}, we used the fact
 that the function $\ell \mapsto \sinh \div{\ell}$ lies in the kernel of the differential operator
 $\cD = \frac 14 - \partial^2$.  Thus, the integral
 $\int_{0}^{+ \infty} \cD H_L(\ell) \sinh \div{\ell} \d \ell$, instead of growing exponentially in
 $L$ as one could expect, is reduced to the bounded term~$-H_L(0)$. Therefore, the integration by
 parts with respect to the operator~$\cD$ has completely solved the problem raised in the previous
 section, namely the exponential contribution of order~$e^{L/2}$ coming from simple geodesics on the
 geometric side and from the trivial eigenvalue on the spectral side.

\begin{summ}
  Applying the operator $\cD$ on the test function $H_L$ both kills the trivial eigenvalue on the
  spectral side and the dominant contribution of simple geodesics in the averaged geometric term.
  The crucial point of the forthcoming analysis is that the kernel of $\cD^m$ consists of functions
  of the form $P(\ell)e^{\ell/2}+ Q(\ell)e^{-\ell/2}$, where $P$ and $Q$ are polynomials of degrees
  $\leq m-1$ (we shall actually only use an inclusion, namely the fact that functions of the form
  $P(\ell)e^{\ell/2}$ are in the kernel of~$\cD^m$).
\end{summ}


\subsection{Volume functions for local types, and the Friedman--Ramanujan property}
\label{s:topo}

Our papers \cite{Ours1, Ours2} show that the cancellation by integration by parts, discovered above,
takes place not only for simple geodesics, but for any local topological type of closed
geodesic. Besides, it occurs not only for the leading order term, but to any order in the asymptotic
expansion in inverse powers of $g$. Let us present these results.

 \subsubsection{Local topological types of loops}
\label{sec:local-topol-types}

 In \cite[\S 4]{Ours1} , we define the notion of {\em{local topological type}} for closed geodesics
 in the surface $X$. This allows to break the sum \eqref{e:def_av_all} into
\begin{align*}
 \label{e:def_av_all}
  \av[\text{all}]{F} = \sum_{\type}   \av[\type]{F}
  \end{align*}
where $\sum_\type$ is a sum over all possible local topological types, and
    \begin{equation*}
  \av[\type]{F} := 
  \Ewp{\sum_{\gamma \sim \type} F(\ell(\gamma))}.
\end{equation*}
We write $\gamma \sim \type$ to mean that $\gamma$ has local topological type $\type$. All simple closed
geodesics form one local topological type -- and the corresponding average is the sum
$\av[\mathrm{s}]{F}$ defined in \eqref{eq:av_S}.

\begin{figure}[h]
  \centering
  \includegraphics[scale=0.4]{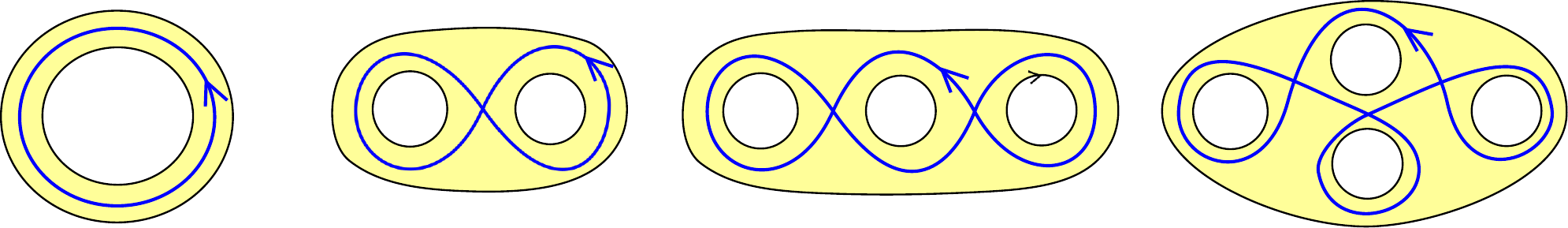}
  \caption{Examples of local topological types. From left to right, we have the type ``simple", a
    figure-eight, and then two generalized eights (see \S \ref{s:gene_eight}).}
  \label{fig:examples_types}
\end{figure}

In a nutshell, a local topological type is an equivalence class $\eqc{\Sf, \curve}$, where $\Sf$ is a
surface (the ``filling type'') and $\curve$ is a loop that fills $\Sf$, for an equivalence
relation that identifies pairs ${{(\Sf, \curve)}}$ having the same topology \cite[\S
4]{Ours1}. Examples of local types are represented in Figure \ref{fig:examples_types}.

  \begin{nota}\label{n:FT}
    To any pair of integers $(g_\Sf, n_\Sf)$ such that $ 2g_\Sf-2+n_\Sf \geq 0$, we shall associate a
    \emph{fixed} smooth oriented surface $\Sf$ of signature $(g_\Sf,n_\Sf)$. We denote by
    $\chi(\Sf)= 2g_\Sf-2+n_\Sf$ the absolute Euler characteristic of $\Sf$.  The data of the pair of
    integers $(g_\Sf, n_\Sf)$, or equivalently of the surface $\Sf$, is called a \emph{filling type}.
  \end{nota}

  We say that a loop $\curve$ fills $\Sf$ if all connected components of $\Sf\setminus \curve$ are either
  contractible or annular regions around a boundary component of $\Sf$ (keep in mind that with this
  convention, we do not allow connected components of $\Sf\setminus \curve$ to be cylinders, unless
  they share a boundary of $\Sf$). The notion extends to what we call \emph{multi-loops} in
  \cite{Ours1, Ours2}, that is, finite collections $\curve=(\curve_1, \ldots, \curve_{K})$ of loops.

  \begin{defa}
  \label{def:local_type}
  A \emph{local loop} is a pair $(\Sf,\curve)$, where $\Sf$ is a filling type and $\curve$ is a
  primitive loop filling $\Sf$. Two local loops $(\Sf,\curve)$ and $(\Sf',\curve')$ are said to be
  \emph{locally equivalent} if $\Sf=\Sf'$ (i.e.  $g_\Sf = g_{\Sf'}$ and $n_\Sf=n_{\Sf'}$), and there
  exists a positive homeomorphism $\psi : \Sf \rightarrow \Sf$, possibly permuting the boundary
  components of $\Sf$, such that $\psi \circ \curve$ is freely homotopic to $\curve'$ (the
  definition implies that the orientation of $\curve'$ is respected). This defines an equivalence
  relation $\eq$ on local loops. Equivalence classes for this relation are denoted as
  $\eqc{\Sf, \curve}$ and called \emph{local (topological) types} of loops. If
  $\type=\eqc{\Sf, \curve}$ is a local topological type of filling type $\Sf$, we let
  $\chi(\type)=\chi(\Sf)$.
\end{defa}

\begin{defa}
  We say that a periodic geodesic $\gamma$ in an arbitrary hyperbolic surface $X$ \emph{has local
    topological type $\type$} if $(\Sf(\gamma), \gamma)$ is a representative of the equivalence class $\type$,
  where $\Sf(\gamma)$ is the surface filled by $\gamma$ in $X$, constructed in \cite[\S 4.3,
  Definition 4.1]{Ours1}. 
\end{defa}

\subsubsection{Volume functions}
\label{sec:volume}

In \cite{Ours1}, we examine the average $\av[\type]{F}$ for a local type $\type$. We prove in
Proposition 5.9 that there exists a unique continuous function
$V_g^\type : \R_{> 0} \rightarrow \R_{\geq 0}$, called \emph{volume function}, such that, for any
measurable function $F$ with compact support,
\begin{equation*}
  \av[\type]{F} = \frac{1}{V_g} \int_{0}^{+\infty} F(\ell) V_g^{\type}(\ell) \d \ell.
\end{equation*}
This density can be computed using our integration formula for the average $\av[\type]{F}$, Theorem 5.5
in \cite{Ours1}, which we now state.
\begin{nota}
  \label{nota:teich}
     For $2-2g-n<0$, we define
   \begin{equation*}
     \mathcal{T}_{g,n}^* := \{(\x, Y), \x \in \R_{>0}^n, Y\in \mathcal{T}_{g ,n}(\mathbf{x})\}
   \end{equation*}
   where $\mathcal{T}_{g ,n}(\mathbf{x})$ is the standard Teichm\"uller space of hyperbolic surfaces
   of genus $g$, with $n$ labelled boundary components of lengths $\x=(x_1, \ldots, x_n)$. This
   space is naturally equipped with the measure $\d \x \d \mathrm{Vol}_{g,n,\x}^{\mathrm{WP}}$,
   where $\d \x$ denotes the Lebesgue measure on $\R_{>0}^n$ and
   $\d \mathrm{Vol}_{g,n,\x}^{\mathrm{WP}}$ the volume induced by the Weil--Petersson form on
   $\cT_{g,n}(\x)$.  We denote by $\mathcal{M}_{g,n}^*$ the corresponding moduli spaces, obtained by
   quotienting the Teichm\"uller space by the action of the Mapping Class Group. Recall that the
   latter is (modulo isotopy) the group of orientation preserving homeomorphisms that fix the
   boundary, so that an element of $\mathcal{M}_{g,n}^*$ still has a labelled boundary.

   If $Z=(\x, Y)\in \mathcal{M}_{g,n}^*$, denote by $\ell(\partial Z)=\sum_{i=1}^n x_i$ the total
   boundary length of $Z$, and $\ell^{\mathrm{max}}(\partial Z)=\max_{i=1}^n x_i$ the length of the
   longest boundary component of $Z$.
\end{nota}
\begin{thm}[\cite{Ours1}, Theorem 5.5]
  \label{thm:density_expr} Let $\type=\eqc{{{\Sf, \curve}}}$ be a local topological type.  For
  any $g \geq 3$, any test function $F$,
  \begin{equation*}
    \av[\type]{F}
    = \frac{1}{m(\type)} \int_{\cT_{g_\Sf,n_\Sf}^*} F(\ell_Y(\curve)) \, \Phi_g^\Sf(\x) \d \x \d
    \mathrm{Vol}_{g,n,\x}^{\mathrm{WP}}(Y)
  \end{equation*}
  where:
  \begin{itemize}
  \item $m(\type)$  is the number of permutations of $\partial \Sf$ stabilising the homotopy class
    of $\curve$;
  \item for any metric $Y$ on the surface $\Sf$, $\ell_Y(\curve)$ denotes the length of the geodesic
    representative in the free homotopy class of $\curve$;
  \item for any $\x = (x_1, \ldots, x_{n_\Sf}) \in \R_{>0}^{n_{\Sf}}$,
    \begin{equation*}
      \Phi_g^\Sf(\x)
      = \frac{x_1 \ldots x_{n_\Sf}}{V_g} \sum_{R \in R_g(\Sf)} V_R(\x)
    \end{equation*}
    with $R_g(\Sf)$ an enumeration of all the possible embedding of $\Sf$ in a surface of genus $g$ (see
    \cite[Definition 4.12]{Ours1}), and $V_R(\x)$ the product of the Weil--Petersson volumes of the
    moduli spaces appearing in the complement of $\Sf$ in this embedding.
  \end{itemize}
\end{thm}
A more explicit example of this formula is given in \S \ref{sec:real-fill-type}, when $\Sf$ is a pair
of pants.

 \subsubsection{Asymptotic expansions in inverse powers of $g$}
 \label{sec:asympt-expans-invers}

 Let us now perform asymptotic expansions, and analyse the average $\av{F}$ at an arbitrary order of
 precision in $g^{-1}$.
 
 The following proposition is useful to give a probabilistic bound on the Euler characteristic of
 the surface filled by a closed geodesic~$\gamma$ on a surface of large genus $g$.

 \begin{prp}\label{p:apriori} The probability for a random hyperbolic
   surface of genus $g$ to contain a multi-loop of length $\leq L$ filling a surface of absolute
   Euler characteristic $> \chi$ is $\O[\chi]{L^{c(\chi)} e^{L} /g^{\chi+1}}$. 
\end{prp} 
Together with Lemma \ref{lem:bound_number_closed_geod}, this allows to decompose the average
$\av[\text{all}]{F} $ into
\begin{align} \label{e:bound_chi} \av[\text{all}]{F} =\sum_{\type,\, \chi(\type) \leq \chic}
  \av[\type]{F}+ \cO_\chic\Big(\frac{L^{c(\chic)} e^{L}}{g^{\chic}} L \norm{F(\ell) e^{\ell}}_\infty\Big)
\end{align}
if $F$ is supported in $[-L, L]$.  If $L=A\log g$, and we choose for instance
$\chic= \lceil 2A \rceil +1$  (without any attempt to be optimal), then the error term is
  $\cO_A( \norm{F}_\infty)$. As a consequence, for our prospects, it is enough to study the term
  $\sum_{\chi(\type) \leq \chic} \av[\type]{F}$ for arbitrarily large fixed $\chic$.

In \cite{Ours1}, we prove the following result, for any individual term $\av[\type]{F}$.
\begin{thm}[\cite{Ours1}, Theorem 5.14]
  \label{thm:exist_asympt_type_intro} Let $\type$ be a local topological type.  The volume function
  $V_g^{\type}$ admits an asymptotic expansion in powers of $g^{-1}$, in the following sense. There
  exists a unique family of continuous functions $(f_k^{\type})_{k \geq \chi(\type) }$ such that,
  for any integer $K \geq 0$, any large enough~$g$, any $\eta>0$,
  \begin{equation*}
    \av[\type]{F} = \sum_{k=\chi(\type)}^K
    \frac{1}{g^k} \int_0^{+ \infty} F(\ell) f_k^{\type}(\ell) \d \ell
    + \O[\Sf,K, \eta]{ \frac{\norm{ F(\ell) \, e^{(1+\eta)\ell}}_\infty}{g^{K+1}}}.
  \end{equation*}
\end{thm}
Note that the expansion starts at the absolute Euler characteristic $\chi(\type) \geq 0$ of $\type$.

As a consequence of this theorem, for any $m\geq 1$, the term
$\av[\type]{\ell \, e^{- \frac \ell 2} \, \mathcal{D}^m H_L(\ell) }$ appearing in Lemma
\ref{lem:selberg_reformulated} may be developed into:
   \begin{equation}\label{e:ouais}
     \sum_{k=\chi(\type)}^K
    \frac{1}{g^k} \int_0^{+ \infty} \ell \, e^{- \frac \ell 2} \,  f_k^{\type}(\ell) \, \mathcal{D}^m
    H_L(\ell) \d \ell
    + \O[\Sf,K, \eta]{ \frac{\norm{ \ell \, e^{(\frac12+\eta)\ell} \, \mathcal{D}^m H_L(\ell)}_\infty}{g^{K+1}}}.
  \end{equation}
  As in the discussion following \eqref{e:bound_chi}, we can pick the order $K$ of the expansion so that its
  remainder is negligible in the context of our proof of \eqref{e:dream}. We therefore are left to
  understand the integrals appearing in \eqref{e:ouais}.
 
\subsubsection{The Friedman--Ramanujan property and integration by parts argument}
\label{sec:friedm-raman-prop}

The central focus of our paper \cite{Ours2} is to understand the structure of the coefficients
$f_k^{\type}(\ell)$: we prove that they satisfy the Friedman--Ramanujan property (in the weak
sense). This property allows to perform an integration by parts in each term of \eqref{e:ouais}.
\begin{defa}\label{def:weakFR}
  Let $c\geq 0$.  A continuous function $f: \R_{\geq 0} \rightarrow \C$ belongs in the set
  $\FRrem_w^c$ if there exists $ c_1>0$ such that, for all $n\geq 1$,
 \begin{align}\label{e:weak}
 \int_0^{n} |f(s)| ds \leq c_1 (n+1)^c e^{\frac n 2}.
 \end{align}
  We denote $\cR_w=\cup_{c\geq 0} \cR_w^c$.
  
  A continuous function $f : \R_{\geq 0} \rightarrow \C$ is said to be a
  \emph{Friedman--Ramanujan function} (in the weak sense), if there exists a polynomial function $p$
  such that
  \begin{equation*}
    f(\ell) - p(\ell) \, e^\ell \in \cR_w.
  \end{equation*}
  The function $\ell\mapsto p(\ell) \, e^\ell$ is called the \emph{principal term} of $f$, and
  $f(\ell) - p(\ell) \, e^\ell $ the \emph{remainder term}. The \emph{degree} of the
  Friedman--Ramanujan function $f$ refers to the degree of the polynomial $p$. We denote by
  $\FR_w^m$ the subset of $\FR_w$ for which $p\in \C_{m-1}[\ell]$, i.e. $p$ has degree (strictly)
  less than $m$.
 \end{defa}
The term \emph{weak} refers to the fact that we express the upper bound \eqref{e:weak} in an
integrated form, and not as a pointwise bound. 
\begin{defa} 
  We define a semi-norm $\| \cdot \|_{\FR_w}$ on $\FR_w$ by setting
  $\| f \|_{\FR_w} := \| p \|_{\ell^\infty}$ for any Friedman--Ramanujan function $f$ of principal
  term $p(\ell) \, e^\ell$, where $\| p \|_{\ell^\infty}$ denotes the maximum modulus of the
  coefficients of the polynomial $p$.
  We define a norm on $\cR_w^c$ as
  \begin{align*}
  \norm{g}_{\cR_w^c} := \sup_{n\geq 1}\frac{\int_0^{n} |g(s)| \d s}{  (n+1)^c e^{\frac n 2}}.
 \end{align*}  
 We let
 $\cF_w^{m, c} := \C_{m-1}[\ell] e^\ell \oplus \cR_w^c$, with norm $\norm{p(\ell)e^{\ell} +g}_{\cF^{m, c}}:= \norm{p}_{\infty}+ \norm{g}_{\cR_w^c}$.  
  \end{defa}
 
  The crucial ingredient to prove \eqref{e:dream} is the following result, which is the main result
  of \cite{Ours2}, and the focus of \S \ref{s:howFR} of this article.
  
  \begin{thm}[{\cite{Ours2}}] \label{t:leviathan} Let $\type=\eqc{{{\Sf, \curve}}}$ be a local
    topological type.  For any $k$, the function $\ell \mapsto \ell^{s_\type} f_k^{\type}(\ell)$ is
    a weak Friedman--Ramanujan function, where $s_\type=1$ if $\type$ corresponds to the local
    type of simple geodesics, and $s_\type=0$ for all other local type.
   
    More precisely, there exists $m(k, \Sf)$ and $c(k, \Sf)$, depending only on the filled surface
    $\Sf$, and not on the loop $\curve$, such that
    $\ell^{s_\type} f_k^{\type}\in \cF_w^{m(k, \Sf), c(k, \Sf)}$, and
   $$\norm{ \ell^{s_\type} f_k^{\type}}_{\cF^{m(k, \Sf), c(k, \Sf)}} \leq c(k, \Sf).$$
   \end{thm}
   
To rephrase Theorem \ref{t:leviathan}, we can write 
  $$\ell^{s_\type} f_k^{\type}(\ell)=p_k^{\type}(\ell) \, e^\ell +  r_k^{\type}(\ell) $$ where
  $p_k^{\type}$ is a polynomial function of degree $d_k^\type\leq m(k, \Sf)-1$, and $r_k^{\type}$
  satisfies \eqref{e:weak} with some constants $c, c_1\leq c(k, \Sf)$.

  Without loss of generality, we can assume that $m(k, \Sf)\leq m(k+1, \Sf)$ and
  $c(k, \Sf)\leq c(k+1, \Sf)$.

  The relevance of the set $\FR_w^m$ of Friedman--Ramanujan functions is that their principal term,
  multiplied by the function $e^{-\ell/2}$ coming from the trace formula, falls in the kernel of the
  power~$\cD^{m}$.  This allows an iterated integration by parts in
  \eqref{e:ouais}, generalizing \eqref{e:vanish}:
  \begin{prp}[{\cite[Proposition 3.18]{Ours1}}]
  \label{cor:FR_implies_small}
  Let $\type=[\Sf, \curve]_{loc}$ be a local topological type. For any integer $K \geq 0$, there
  exists constants $c(K, \Sf), m(K, \Sf) \geq 0$ (depending only on the filled surface~$\Sf$ and not
  the loop $\curve$) such
  that, for any large enough $g$, any $m \geq m(K, \Sf)$, any $\eta >0$ and any $L \geq 1$,
  \begin{equation}\label{e:contribTIBP}
    \avb[\type]{\ell \, e^{- \frac \ell 2} \, \D^m H_L(\ell)}
    = \O[m,\Sf,K,\eta]{L^{c(K, \Sf)} + \frac{e^{L(\frac{1}2 +\eta)}}{g^{K+1}}}.
  \end{equation}
\end{prp}

As a result, the contribution of a given local topological type to the quantity
$\avb{\ell \, e^{- \frac \ell 2} \, \mathcal{D}^m H_L(\ell) }$ in Lemma
\ref{lem:selberg_reformulated} is much smaller than expected: it is polynomial in $L$, rather than
exponential.  Note that to balance the last error term, of size $\frac{e^{L/2}}{g^{K+1}}$, with our
other error terms $\cO(g)$ coming from the topological term in the trace formula, we are lead to
take $L=2(K+2)\log g$. With this choice, \eqref{e:contribTIBP} reads
  \begin{equation*}
    \avb[\type]{\ell \, e^{- \frac \ell 2} \, \D^m H_L(\ell)}
    = \O[m,\Sf,K,\eta]{ g^{1+2\eta(K+2)}} = o(e^{(\alpha+\eps)L})
  \end{equation*}
  with $\alpha =\frac1{2(K+2)}$ and $\eta <\eps$, which is exactly the bound we need to prove \eqref{e:dream}.

  \begin{summ}
    So far, our approach to prove that $\Pwp{\lambda_1 \leq \frac{1}{4} - \alpha^2 - \epsilon} \To 0$
    consists in the following steps.
    \begin{enumerate}
    \item Using Markov's inequality, the Selberg trace formula and our cancellation argument for the
      trivial eigenvalue, we reduce the problem to finding a scale $L$ and an index $m$ such that
      the average $\av{\ell e^{-\ell/2}\D^m H_L(\ell)}$ is $o(e^{(\alpha + \epsilon)L})$ as
      $g \rightarrow + \infty$ (see Lemma \ref{lem:selberg_reformulated}).
    \item We compute the asymptotic expansion for the average $\av{\ell e^{-\ell/2} \D^mH_L}$ in
      inverse powers of $g$ at the appropriate order $K$, so that $\frac1{2(K+2)}\leq \alpha$.
    \item We pick the length-scale $L=A\log g$ for $A=2(K+2)$ so that remainders of the asymptotic
      expansion $e^{L/2}/g^{K+1}$ are of the same size as the inevitable topological term of size
      $g$.
    \item We reduce ourselves to studying only the filling types $\Sf$ of absolute Euler
      characteristic $\chi(\Sf) \leq \chic$ for $\chic=2A+1$ by using equation \eqref{e:bound_chi}.
    \item \label{i:ibp_term} Thanks to the integration by part argument in Proposition
      \ref{cor:FR_implies_small}, we can find an integer $m$ such that, for each type $\type$ such that
      $\chi(\type) \leq \chic$, the individual term $\av[\type]{\ell e^{-\ell/2}\D^m H_L(\ell)}$ satisfies the
      desired bound.
    \end{enumerate}
    Note that the case $K=0$ corresponds to $\alpha = \frac 16$, treated in
    \cite{wu2022,lipnowski2021}, and the case $K=1$ to our work \cite{Ours1}. The optimal spectral
    gap, \eqref{e:dream}, requires $K$ to be arbitrarily large.
  \end{summ}
  However, this is not enough to conclude to the proof of \eqref{e:dream}, because whilst we can
  bound each individual term in step \eqref{i:ibp_term}, the number of types $\type$ appearing in the
  average $\av{\ell e^{-\ell/2}\D^m H_L(\ell)}$ is itself exponential in $L$. As a consequence, it
  is impossible to sum the bounds from \eqref{i:ibp_term}, and another argument needs to be
  found. We present this challenge in the next subsection, and its resolution in \S \ref{s:tangles_moebius}.

   \subsection{Tangles and how to remove them}
   \label{s:tangle_pb}

   \subsubsection{The failure of the trace method}
   \label{sec:necess-remov-tangl}

   In \cite[Theorem 9.1]{Ours1}, we prove the following negative result about the sum
   $f_1^{\mathrm{all}} := \sum_\type f_1^\type$.
   \begin{prp}
     \label{p:FR_not}
     The function $\ell \mapsto \ell f_1^{\mathrm{all}}(\ell)$ is not a weak Friedman--Ramanujan
     function.
   \end{prp}
   So, whilst Theorem \ref{t:leviathan} tells us that each individual function
   $\ell \mapsto \ell f_1^\type(\ell)$ is a Friedman--Ramanujan function, the sum over all local types
   $\type$ is not. This might seem surprising because of the stability of the Friedman--Ramanujan
   property by finite addition. However, the set of local types $\type$ is countable, and (without any
   further hypothesis) the number of local types appearing when examining geodesics of length
   $\leq L$ is exponential in~$L$. This exponential proliferation of local types is the reason why
   the Friedman--Ramanujan property does not hold once we sum over all local types.
   
   Unfortunately, this means that the trace method sketched so far is doomed to fail and cannot be
   used \emph{as is} to conclude to the proof of \eqref{e:dream}.  We can only apply our integration
   by parts argument, Proposition~\ref{cor:FR_implies_small}, to individual terms corresponding to a
   fixed local type $\type$, but not to the sum over all local types (which is what appears in the
   Selberg trace formula).

   The proof of Proposition~\ref{p:FR_not} consists in putting together the following two
   observations, Lemmas 9.4 to 9.6 in \cite{Ours1}.
   \begin{itemize}
   \item There exists a constant $c>0$ such that, for any small enough $a>0$,
     \begin{equation}
       \label{eq:P_small_eig}
       \Pwp{\lambda_1 \leq a} \geq c \, \frac{a}{g} \cdot
     \end{equation}
   \item If $\ell \mapsto \ell f_1^{\mathrm{all}}(\ell)$ were a Friedman--Ramanujan function, then
     we could use the trace method to prove that
     \begin{equation}
       \label{eq:P_if_FR}
       \Pwp{\delta \leq \lambda_1 \leq \frac{5}{72}} = \O[\delta]{g^{-\frac 54}}
     \end{equation}
     which is false since $1/g \gg 1/g^{5/4}$.
   \end{itemize}
   Note that the failure of the trace method is related to the fact that the probability for the
   spectral gap $\lambda_1$ to be small is \emph{not small enough}: it is an event of probability
   going to $0$ as $g \rightarrow + \infty$, but only at a speed $1/g$.
   
   \subsubsection{Tangles and the exponential proliferation of local types}
\label{sec:role-tangl-expon}

Equation \eqref{eq:P_small_eig} is proven by exhibiting a set of surfaces for which $\lambda_1 \leq
a$: surfaces which contain an embedded once-holed torus of boundary length $\leq a$. It is easy to
use Mirzakhani's methods to estimate the probability of this geometric event.

We refer to embedded once-holed tori with a short boundary as \emph{tangles}: they are geometric
patterns responsible for difficulties in our trace method, or more generally in the study of the
spectral gap $\lambda_1$. The other patterns that we refer to as tangles are embedded pairs of pants
with short boundary as well as extremely short closed geodesics -- see \S \ref{s:deftangles} for a
definition. The denomination ``tangles'' is taken from the work of Friedman, and later Bordenave,
who defined tangles in regular graphs \cite{friedman2003, bordenave2020}.  Our definition, close to
the geometric definition of Bordenave, has been developed and investigated by the second author and
Thomas in \cite{monk2021a}, and appears in Lipnowski--Wright's paper \cite{lipnowski2021}.

\emph{Tangle-free surfaces} are surfaces that contain no tangles; we shall denote this set as $\tf$
in this subsection for the sake of simplicity (in reality, this set depends on parameters $\kappa$
and $\tfL$ related to our parameter $\alpha$, and will be denoted as $\Atf$). We make the following
observations:
\begin{itemize}
\item $\Pwp{\tf}$ goes to $1$ as $g \rightarrow + \infty$, so that we can safely condition the
  Weil--Petersson probability model by the tangle-free hypothesis to prove \eqref{e:dream};
\item for any $X \in \tf$, we can prove polynomial bounds on the number of local types appearing in
  the trace method (see Theorem \ref{cor:TF_curves}).
\end{itemize}
As a consequence, the tangle-free hypothesis is a good solution to remove the issue of exponential
proliferation of local types. This conditioning can easily be performed by multiplying the trace
formula by the indicator function $\1{\tf}$ in the trace method, in which case the quantity that we
need to describe (from Lemma \ref{lem:selberg_reformulated}) becomes
\begin{equation}
  \label{eq:trace_after_TF}
  \Ewp{\sum_{\gamma \in \geod(X)} \frac{\ell(\gamma)
      \D^m H_L(\ell(\gamma))}{2 \sinh \div{\ell(\gamma)}} \1{\tf}(X)}.
\end{equation}
Importantly, the number of local types in this sum is at most polynomial in $L$, which will allow us to
study the contribution of local types individually. We will therefore aim to extend the results of
\S \ref{s:topo}, and in particular Theorem \ref{t:leviathan}, to the density functions $V_g^{\type,\tf}$
defined by
\begin{equation}
  \label{e:density_vol_TF}
  \int_0^{+\infty} F(\ell) V_g^{\type,\tf}(\ell) \d \ell
  = \Ewp{\sum_{\gamma \sim \type} F(\ell(\gamma)) \1{\tf}(X)}.
\end{equation}
   
\subsubsection{Inclusion-exclusion and Moebius formula}
\label{sec:incl-excl-moeb}

Unfortunately, Mirzakhani's integration formula only allows to compute certain types of averages:
expectations of geometric counting functions. They cannot directly be used to study the
average~\eqref{e:density_vol_TF}.
\begin{rem}
  \label{rem:positivity}
  Had we been averaging functions of constant sign, we could have hoped to use the trivial bound
  $\1{\tf} \leq 1$ to remove the indicator function at this stage. However, the discussion in \S
  \ref{s:killing} shows how the trivial eigenvalue forces us to use non-trivial cancellations and
  oscillatory functions, which makes this impossible.
\end{rem}

A workaround that has been used in \cite{mirzakhani2013,lipnowski2021} consists in performing an
inclusion-exclusion, hence rewriting the indicator function $\1{\tf}$ as a combination of geometric
counting functions. The idea is the following.  Let $N$ be a random variable counting the number of
appearances of a geometric pattern, e.g. short closed geodesics. Then, it is classic that we can
write
\begin{equation}
  \label{e:incl_excl}
  \1{N=0} = 1 - \sum_{j=1}^{+ \infty} \frac{(-1)^j}{j!} N_j
\end{equation}
where, for any $j$, $N_j$ denotes the number of ordered families of $j$ patterns counted by
$N$. This is a rewriting of the indicator function $\1{N=0}$ as a sum of counting functions, which
is what we wished to apply Mirzakhani's integration formula.

In the case where $N$ counts closed geodesics of length $\leq \epsilon$, like in
\cite{mirzakhani2013,lipnowski2021}, thanks to the collar lemma, we can take $\epsilon$ small enough
so that $N_j$ counts \emph{multicurves}, i.e. families of $j$ disjoint simple geodesics of lengths
$\leq \epsilon$. This greatly simplifies the enumeration of all possible topologies for the $j$
geodesics in $N_j$. However, for our purposes, we need $N$ to count tangles, which include embedded
tori and pairs of pants. This is significantly more complicated, because we a priori need to
enumerate all possible topologies of families of $j$ tangles. This enumeration is done explicitly
in \cite{Ours2} at the precision $1/g^2$, where it is still manageable; however, in general, it is
unrealistic to perform it fully.

In \S \ref{s:moebius}, we present a Moebius formula, Corollary \ref{p:moebius-goal}. It is a
rewriting of the indicator function $\1{\tf}$ similar to \eqref{e:incl_excl} but less explicit,
which hence allows to estimate the average~\eqref{e:density_vol_TF} without any topological
enumeration. We then explain how to conclude in \S \ref{s:obtention}. Note that \S \ref{s:howFR} and
\S \ref{s:tangles_moebius}-\ref{s:obtention} are independent, and that the interested reader can
skip \S\ref{s:howFR} to obtain the full walkthrough of our trace method.

\begin{summ}
  The issue of exponential proliferation of local types is solved by making the tangle-free
  assumption. This causes a new challenge: we now need to compute averages containing the indicator
  function of the set of tangle-free surfaces. The Moebius formula developed in \cite{Moebius-paper}
  allows us to rewrite this indicator function in a manageable form, and to conclude the proof
  of~\eqref{e:dream}.
\end{summ}

\section{How to prove the Friedman--Ramanujan property } 
\label{s:howFR}

In this section, we dive in more detail in the proof of Theorem \ref{t:leviathan} and present its
main steps.

\subsection{Friedman--Ramanujan functions as solutions of an integro-differential equation}

We define two operators acting on locally integrable functions, $\primitive=\primitive_x=\int_0^x$
(taking the primitive vanishing at $0$) and $\cL=\cL_x= \Id -\int_x=\Id-\primitive_x$ (where $\Id$
stands for the identity operator). It is trivial, nevertheless useful, to see that $\primitive$ (and
hence also $\cL$) preserves the classes of functions $\cF, \cF_{w}, \cR, \cR_w$. More precisely,
$\cL$ is a continuous operator from $\cF^m$ to $\cF^{m-1}$, and from $\cR^c$ to itself (and
same goes with the weak spaces).
  
A fundamental ingredient in the proof of Theorem \ref{t:leviathan} is the fact that Friedman--Ramanujan functions are characterized as solutions (modulo $\cR$) of a certain integral equation:
 \begin{prp}[{\cite{Ours2}}] \label{p:charFR} \quad 
\begin{itemize}
\item $f\in \FR^m$ if and only if $\cL^{m} f\in \cR$.
\item$f\in \FR^m_w$ if and only if $\cL^{m} f\in \cR_w$.
\end{itemize}

\end{prp}

\subsection{Stability under convolution and almost-convolution}
A second important fact behind the proof of Theorem \ref{t:leviathan} is the stability of the
classes $\FR, \FR_w,\cR, \cR_w$ under convolution.  This was proven in \cite[Proposition 3.5]{Ours1}
by direct calculation. Another possible approach to understand this stability is via the usual property
of convolution (well known for differential operators, but valid also for the operator $\cL$): given
two functions $f_1, f_2$ and integers $m_1, m_2$, we have
\begin{align}\label{e:easyconvo}
\cL^{m_1+m_2} (f_1 * f_2)= \cL^{m_1} f_1 * \cL^{m_2} f_2.
\end{align}
If $f_i\in \cF^{m_i}$, then $\cL^{m_i} f_i \in \cR$ for $i=1, 2$. Thus, using the stability
of~$\cR$ under convolution (a straightforward upper bound), we obtain that
$\cL^{m_1+m_2} (f_1 * f_2)\in \cR$, i.e.  $f_1 * f_2\in \cF^{m_1+m_2}$.

We shall use the following extension of the notion of convolution:
 \begin{defa}[$(h, \varphi)$-convolution, \cite{Ours2}] 
   Let $\varphi :\IR^n\To \IR$ be a measurable, locally bounded function, and let $h$ be a
   measurable function defined on the support of $\varphi$. We assume that for every compact subset
   $K\subset \IR$, the set $h^{-1}(K) \cap \{\varphi \not = 0\}$ is bounded.

   We call $(h, \varphi)$-convolution of $n$ locally integrable functions
   $f_1, \ldots, f_n: \IR \To \IR$ the pushforward of the distribution
   $\varphi(x_1, \ldots, x_n)\prod_{i=1}^n f_i(x_i) \d x_i$ by the map $h$. In other words, it is
   the distribution $\mu$ on $\IR$ defined by
 \begin{align*}
   \int_{ \IR} F(x) \, \mu(\d x)
   = \int_{\IR^n}F(h(x_1, \ldots, x_n)) \, \varphi(x_1, \ldots, x_n) \prod_{i=1}^n f_i(x_i) \d x_i
 \end{align*}
 for every continuous compactly supported function $F$.
 \end{defa}

 \begin{nota}
   Under relatively mild smoothness conditions on the function $h$, the distribution $\mu$ admits a
   density, which will be denoted $f_1\star \ldots \star f_n |^h_\varphi$.
 \end{nota}

The usual convolution corresponds to $h( x_1, \ldots,  x_n)= x_1+\ldots+ x_n$, and either $\varphi\equiv 1$, or $\varphi( x_1, \ldots , x_n)=\bbbone_{ x_1>0, \ldots,  x_n>0}$ (in the case of convolutions
of functions defined on $\R_{\geq 0}$).
 In what follows, we fix $\Amin>0$ and we will be interested in functions defined on $[\Amin, +\infty)^n$.
We define a space $E_n^{(\Amin)}$ of functions which are ``exponentially close to being constant'', and an affine space $\cE_n^{(\Amin)}$ of functions
which are ``exponentially close to the function $L_n(x_1,\ldots,x_n)=x_1+\ldots+x_n$''.

\begin{defa}[\cite{Ours2}]
  We let $E_n^{(\Amin)}$ be the vector space of functions $h: [\Amin, +\infty)^n \rightarrow \R$
  such that
  \begin{itemize}
  \item $h$ is of class $\cC^\infty$;
  \item for every multi-index $\alpha\in \{0, 1\}^n$ such that $\alpha\not=(0,\ldots, 0)$, if
    $\alpha \cdot \x :=\sum_{i=1}^n \alpha_i x_i $,
    \begin{align*}
      \norm{h}_{\alpha, \Amin} :=
      \sup_{\x \in [\Amin, +\infty)^n}\Big\{ e^{\alpha \cdot \x} \,|\partial^\alpha h(\x)|\Big\} <+\infty.
    \end{align*}
  \end{itemize}
  We let $\cE_n^{(\Amin)}$ be the affine space of functions $h$ such that $h -L_n \in E_n^{(\Amin)}$.
  Equivalently,
\begin{itemize}
\item $h$ is of class $\cC^\infty$;
\item for every multi-index $\alpha\in \{0, 1\}^n$ such that $|\alpha| = \alpha_1 + \ldots + \alpha_n \geq 2$, 
\begin{align*} 
\sup_{[\Amin, +\infty)^n}\Big\{ e^{\alpha \cdot \x}|\partial^\alpha h(\x)|  \Big\} <+\infty;
\end{align*}
\item for every $j\in \{1, \ldots, n\}$,  
\begin{align*} 
\sup_{[\Amin, +\infty)^n}\Big\{ e^{  x_j}|\partial_j h(\x)-1| \Big\} <+\infty.
\end{align*}
\end{itemize}
 \end{defa}
 In geometric applications, our main examples of such functions will be the following.
  \begin{nota}
   For $x\in \R$, we denote
   $$\hyp_{-}(x)=\sinh(x), \quad \hyp_{+}(x)=\cosh(x), \quad\hyp_{0}(x)=1.$$
   For $\x=(x_1,\ldots, x_n)\in \R^n$ and $\eps=(\eps_1, \ldots, \eps_n)\in \{-, 0, +\}^n$, we write
   \begin{align}
     \label{e:hypnotation}
     \hyp_\eps(\x) := \prod_{i=1}^n \hyp_{\eps_i}(x_i).
   \end{align}
 \end{nota}
 \begin{nota}
   \label{nota:hL}
      For $\beta=(\beta_\eps)_{\eps\in \{\pm\}^n}$ a family of non-negative real numbers such that
   \mbox{$\beta_{(1, \ldots, 1)}=1$} and $\x=(x_1,\ldots, x_n)\in \R^n$, we define
 $$h_\beta(\x) := 2\argcosh\Big(\sum_{\eps\in \{\pm\}^n} \beta_\eps \hyp_\eps\div{\x}\Big).$$
\end{nota}
 \begin{prp}\label{p:mainclass}
  For any $\Amin>0$,  $h_\beta \in \cE_n^{(\Amin)}$.  
 Moreover, there exists $C=C(\Amin, n)$ such that 
 \begin{align}\label{e:defEM}
\sup_{\beta}  \sup_{\alpha\not= (0, \ldots, 0)}\norm{h_\beta-L_n}_{\alpha, \Amin} \leq C.
 \end{align}
 \end{prp}

 We prove the stability of the class of Friedman--Ramanujan functions by almost-convolution, which
 will be a key ingredient to our approach to \eqref{e:dream}.
 
 \begin{thm}[{\cite{Ours2}}] \label{t:intermediate} Let $\Amin>0$. We assume that
   $\supp \varphi\subset [\Amin, +\infty)^n$, $h\in \cE_n^{(\Amin)}$, and $\varphi\in E_n^{(\Amin)}$. If for
   $1 \leq i \leq n$, $f_i \in \FR^{m_i}$, then
   $ f_1\star \ldots \star f_n|^h_\varphi\in \FR^{m_1+\ldots+m_n}$.
\end{thm}

The proof consists in using the characterization of Friedman--Ramanujan functions by an integral
equation, as above for the simple case of convolutions. More precisely, we show that
$\cL^{m_1+\ldots +m_n}(f_1\star \ldots \star f_n|^h_\varphi)$ coincides with
$\cL^{m_1}f_1\star \ldots \star \cL^{m_n} f_n|^h_\varphi$ modulo $\cR$, and that the modified
convolution $\cL^{m_1}f_1\star \ldots \star \cL^{m_n} f_n|^h_\varphi$ itself belongs to $\cR$.

Theorem \ref{t:intermediate} is also true in the weak sense, replacing the space $\FR$ by $\FR_w$.

\subsection{Geometric application: the case of the figure-eight}

Let us detail our method to prove Theorem \ref{t:leviathan} in one simple case, the figure-eight
filling a pair of pants. This situation has already been solved very explicitly in \cite[\S
7]{Ours1}; however, here, we use a more general method, which can be extended to all cases.

\subsubsection{Expression of the volume function for loops filling a pair of pant}
\label{sec:real-fill-type}
 
Let $\curve$ be a loop filling a pair of pants $\mathbf{P}$, i.e. a fixed surface of signature
$(0,3)$ with boundary components labelled $1, 2, 3$. Then, the length of $\curve$ is an analytic
function of the lengths of the three boundary components $(x_1, x_2, x_3)$ of $\mathbf{P}$, which we shall
denote as $\mathfrak{L}_\curve : \R_{>0}^3 \rightarrow \R_{> 0}$. In \cite[Example 5.7]{Ours1}, we
explicit the integration formula from Theorem \ref{thm:density_expr}, which yields:
  \begin{equation*}
    \av[\eqc{\mathbf{P}, \curve}]{F}
    = \frac{1}{m(\curve)}
    \iiint_{\R_{>0}^3} F(\mathfrak{L}_\curve(x_1, x_2, x_3)) \, \Phi_g^{\mathbf{P}}(x_1, x_2, x_3)
    \d x_1 \d x_2 \d x_3
  \end{equation*}
  where $m(\curve) \in \{1, 2, 3, 6\}$ is the number of permutations of $\{1, 2, 3\}$ stabilising
  the homotopy class of $\curve$, and 
  \begin{equation} \label{e:realization}
    \begin{split}
      \Phi_g^{\mathbf{P}}(\x)
      =  \frac{x_1 x_2 x_3}{V_g}
      & \left[ V_{g-2,3}(x_1,x_2,x_3)
      + \sum_{g_1+g_2+g_3=g} V_{g_1,1}(x_1) V_{g_2,1}(x_2) V_{g_3,1}(x_3)
      \right.  \\
      & \left.
       + \sum_{\substack{\{i_1, i_2, i_3\} \\ = \{1,2,3\}}}
      \left( \frac{\delta(x_{i_1} - x_{i_2})}{x_{i_1}} V_{g-1,1}(x_{i_3})
      + \sum_{i=1}^{g-2}   V_{i,2}(x_{i_1}, x_{i_2}) V_{g-i-1,1}(x_{i_3}) \right)
      \right].
    \end{split}
  \end{equation}
  The various terms in \eqref{e:realization} correspond to all mapping class group orbits of a pair of pants
  inside a surface of genus $g$. In particular,
  \begin{itemize}
  \item the term $ V_{g-2,3}(x_1,x_2,x_3)$ corresponds to a pair of pants whose complement is a
    connected surface of genus $g-2$ with $3$ boundary components;
  \item $\delta$ denotes the Dirac delta distribution, and the term
    $\frac{\delta(x_{i_1} - x_{i_2})}{x_{i_1}} V_{g-1,1}(x_{i_3})$ corresponds to a once-holed
    torus, whose complement has genus $g-1$ and $1$ boundary component;
  \item all other terms enumerate pairs of pants that disconnect the surface, into two or three
    connected components.
  \end{itemize} 
This shows that the volume function $V_g^{\eqc{\mathbf{P}, \curve}}$ is defined as the unique function such that
\begin{align} \label{e:explicitV}
\int_0^{+\infty} F(\ell)V_g^{\eqc{\mathbf{P}, \curve}}(\ell)\d \ell= 
 \frac{1}{m(\curve)}
    \iiint_{\R_{>0}^3} F(\mathfrak{L}_\curve(x_1, x_2, x_3)) \, \Phi_g^{\mathbf{P}}(x_1, x_2, x_3)
    \d x_1 \d x_2 \d x_3
\end{align}
for any test function $F$. In theory, if we know $\mathfrak{L}_\curve$ and all the volume polynomials explicitly, we can obtain the expression of $V_g^{\eqc{\mathbf{P}, \curve}}$.

\subsubsection{Asymptotic expansion of volume polynomials}

The asymptotic expansion in Theorem \ref{thm:exist_asympt_type_intro} follows from the
asymptotic expansion on the function $\Phi_g^{\mathbf{P}}$ given below (which admits a generalization to any
filling type $\Sf$, \cite[Proposition 5.21]{Ours1}).

\begin{prp}[{\cite[Corollary 1.4]{anantharaman2022}}]
  \label{lem:claim_exp_phi_T}
  There exists a unique family of distributions $(\psi_k^{\mathbf{P}})_{k \geq 1}$ such that, for
  any $K \geq 0$, any large enough~$g$, any $\x \in \R_{>0}^{3}$,
  \begin{equation}
    \label{e:claim_exp_phi_T}
    \Phi_g^{\mathbf{P}}(\x) = \sum_{k=1}^K \frac{\psi_k^{\mathbf{P}}(\x)}{g^k}
    + \Ow[K]{\frac{(\norm{\x}+1)^{c_K^{\mathbf{P}}}}{g^{K+1}} \exp{\div{x_1+x_2+x_3}}}.
  \end{equation}
  for a constant $c_K^{\mathbf{P}}$. 
  Furthermore, for all $k$, the distribution $\psi_k^{\mathbf{P}}$ is odd in each variable and can be uniquely written as a
  linear combination of distributions of the following two forms:
  \begin{align}
    \label{eq:form_coeff_phi}
     & x_1^{k_1} x_2^{k_2} x_3^{k_3}
    \hyp_{\eps} \div{\x} \\
    \label{eq:form_coeff_phi_delta}
     & x_{i_2} \delta(x_{i_1} - x_{i_2})
   \, x_{i_3}^{k_3} \,\hyp_{\eps_3}\div{x_{i_3}},  
  \end{align}
  where $\{i_1, i_2, i_3\}=\{1, 2, 3\}$ and $\eps_i\in\{-, 0, +\}$.
  Furthermore, the leading order term is
   \begin{equation}\label{e:psi_1}
\psi_1^{\mathbf{P}}(\x)=  \sinh \div{x_1}   \sinh \div{x_2}  \sinh \div{x_3} + \sum_{\substack{\{i_1, i_2, i_3\} \\ = \{1,2,3\}}}  x_{i_2} \delta(x_{i_1} - x_{i_2}) \sinh \div{x_{i_3}}.
 \end{equation}
  \end{prp}
  
  The asymptotic expansion of $V_g^{\eqc{\mathbf{P}, \curve}}$ comes from inserting the asymptotic
  expansion \eqref{e:claim_exp_phi_T} into \eqref{e:explicitV}. For instance, by \eqref{e:psi_1},
  the first term is characterized by the relation
  \begin{align}\nonumber
\int_0^{+\infty} F(\ell) f_1^{\eqc{\mathbf{P}, \curve}}(\ell)\d \ell  \hspace{4em} &\\
=\label{e:main} \iiint_{\R_{>0}^3} F(\mathfrak{L}_\curve(x_1, x_2, x_3)) \,  
\Big[ &\sinh \div{x_1}   \sinh \div{x_2}  \sinh \div{x_3} \\
 \label{e:handle}& + \sum_{\substack{\{i_1, i_2, i_3\} \\ = \{1,2,3\}}}  x_{i_2} \delta(x_{i_1} - x_{i_2}) \sinh \div{x_{i_3}}\Big]
    \d x_1 \d x_2 \d x_3
  \end{align}
for any test function $F$.

In what follows, we consider the simplest loop filling a pair of pants, that is to say, the figure-eight, going around the boundary curves numbered $1$ and $2$ of $\mathbf{P}$. In this case,
\begin{align}\label{e:L8}\mathfrak{L}_8(x_1, x_2, x_3)= 2\argcosh\Big(2\cosh\div{x_1} \cosh\div{x_2}+ \cosh\div{x_3}\Big).
\end{align}
We explain how to prove that  $f_1^{\eqc{\mathbf{P}, \curve}}$ is a Friedman--Ramanujan function.

It can be noted that, from the variational characterization of geodesics, or from the formula \eqref{e:L8}, 
$$ x_3\leq \mathfrak{L}_8(x_1, x_2, x_3) \quad \text{and} \quad  x_1+x_2\leq \mathfrak{L}_8(x_1, x_2, x_3) .$$
As a result, $\sinh \div{x_{i_3}}\leq e^{\mathfrak{L}_8/2}$. Therefore, the contribution of the term
\eqref{e:handle} to $f_1^{\eqc{\mathbf{P}, \curve}}$ (which came from the figure-eight in a once-holed
torus) can directly be seen to be in $\cR_w$.  Thus, it suffices to focus on the contribution of the
term \eqref{e:main}.

\subsubsection{A nice change of variable}
We are left with showing that the function $f$, defined by
\begin{equation} \label{e:main8}
  \begin{split}
    \int_0^{+\infty} F(\ell) f(\ell)\d \ell
    = \iiint_{\R_{>0}^3} F(\mathfrak{L}_8(x_1, x_2, x_3)) \,  
                           \sinh \div{x_1}   \sinh \div{x_2}  \sinh \div{x_3} \d \x
  \end{split}
\end{equation}
for any test function $F$, is a Friedman--Ramanujan function. This case is sufficiently simple to
allow for the function $f$ to be computed explicitly, but instead, we illustrate our general method
on this example.

Let us introduce a new variable $T$, defined as the length of the orthogeodesic between the boundary
curves numbered $1$ and $2$ of the pair of pants $\mathbf{P}$.  Elementary hyperbolic trigonometric yields
the well-known formulas:
\begin{align}
  & \cosh(T)=\frac{\cosh\div{x_1}\cosh\div{x_2} +\cosh\div{x_3} }{\sinh\div{x_1}\sinh\div{x_2}} \nonumber\\
  & \label{e:eightformula}\mathfrak{L}_8(x_1, x_2, x_3)= 2\argcosh\Big(\cosh\div{x_1} \cosh\div{x_2}+\cosh(T) \sinh\div{x_1}\sinh\div{x_2} \Big).
\end{align}
If we perform the change of variable $(x_1, x_2, x_3)\mapsto (T, x_1, x_2)$, 
the determinant of the Jacobian admits the simple expression
\begin{equation}
  \label{e:magiccv}
  \sinh \div{x_1}   \sinh \div{x_2} \sinh \div{x_3}   \d x_1 \d x_2 \d x_3=2   \sinh^2 \div{x_1}   \sinh^2 \div{x_2}  \sinh T \d T  \d x_1 \d x_2.
\end{equation}

As a consequence, if we define
  \begin{equation*}
    h_T(x_1, x_2) := 2\argcosh\Big(\cosh\div{x_1} \cosh\div{x_2}+\cosh(T) \sinh\div{x_1}\sinh\div{x_2}
    \Big),
  \end{equation*}
  then for any test function $F$ we have
  \begin{equation} \label{e:cv}
    \int_0^{+\infty} F(\ell) f(\ell)\d \ell
    =  2 \iiint_\domain F (h_T(x_1, x_2)) \,  
    \sinh^2 \div{x_1}   \sinh^2 \div{x_2}   \d x_1 \d x_2  \sinh T \d T
 \end{equation}
  where the domain $\domain$ is
 the range of the variables $(T, x_1, x_2)$. In other words, this means that $f$ can be expressed
 thanks to our extended notion of convolution:
 \begin{lem}  \label{e:Vconvol}
   We have
   \begin{align} 
     f= 2 \int_{T=0}^{+\infty}  f_1\star f_2|^{h_T}_{\varphi_T} \sinh T \d T
   \end{align}
   where $f_1(x)=f_2(x)=\sinh^2 \div{x}$, and $\varphi_T(x_1, x_2)=\bbbone_{\domain}(T, x_1, x_2)$.
\end{lem}
Let us check the hypotheses of Theorem \ref{t:intermediate}, in order to prove that $f \in \FR_w$.
\begin{itemize}
\item Obviously, $f_1$ and $f_2$ are Friedman--Ramanujan functions of degree $0$.
\item For any $T$ and
  any $\Amin>0$, the function $h_T$ belongs to $\cE_2^{(\Amin)}$, thanks to Proposition
  \ref{p:mainclass}.
\item The function $\varphi_T$ is not exactly in $E_2^{(\Amin)}$. However, it is
  constant on $\domain$, and one can check that its singularity at the boundary of $\domain$
  corresponds to the locus where $x_3=0$, and yields a contribution which lies in~$\cR_w$.
\end{itemize}
We cannot apply directly Theorem \ref{t:intermediate} because of the presence of the extra parameter
$T$. However, its method of proof allows to prove that $\cL^{2} f$ differs from
 \begin{align} \label{e:convolT}
\int_{T=0}^{+\infty} \cL f_1\star \cL f_2|^{h_T}_{\varphi_T} \, \sinh T \d T
 \end{align}
 by an error belonging to $\cR_w$. We know that $\cL f_1$ and $\cL f_2$ both belong to $\cR_w$: this is
 enough to show that the function defined by \eqref{e:convolT} is in $\cR_w$. We can then conclude
 that $f$ lies in $\FR_w$. 

 The following two remarks are more technical highlights of some aspects of the general
 case, and can be skipped at first read.
 
 \begin{rem}
   \label{rem:neutral_eight}
   The role of the variables $x_1$, $x_2$ and $T$ are different in the discussion above: while we
   need to perform our integration argument with respect to the two variables $x_1$ and $x_2$ (by
   applying the operator $\cL$ on $f_1$ and $f_2$), the variable $T$ is an untouched integration
   parameter. In~\cite{Ours2}, we call $T$ a \emph{neutral parameter}. The reason for this
   difference is that $T$ ``counts twice'' in the length formula \eqref{e:eightformula}.
   Geometrically, this stems from the fact that the figure-eight is homotopic to a loop going once
   around each boundary geodesic of respective length $x_1$ and $x_2$, and twice along the
   orthogeodesic of length $T$.
 \end{rem}

 \begin{rem}
   \label{rem:boundary_comp}
   In the case of the figure-eight, the above proof can be generalized rather straightforwardly to
   prove that the functions $f_k^{\eqc{\mathbf{P},\curve}}$ belong to $\FR$ for every $k \geq
   1$. Unfortunately, as soon as $k \geq 2$, additional polynomial factors
   $x_1^{k_1}x_2^{k_2}x_3^{k_3}$ appear in \eqref{eq:form_coeff_phi}, which were absent from
   \eqref{e:psi_1} when $k=1$. The contributions $x_1^{k_1}$ and $x_2^{k_2}$ can be viewed as part
   of the Friedman--Ramanujan functions $f_1$ and $f_2$ and cause no further difficulties.  However,
   the factor $x_3^{k_3}$ needs to be examined carefully; this is done by writing $x_3$ as a
   function of the integration variables $x_1$, $x_2$ and $T$, and proving it has small derivatives
   w.r.t. $x_1$ and $x_2$ (so that it is close to a polynomial in $x_1$ and $x_2$).  Precise
   estimates on $x_3$ can be found in \cite[Lemma 7.11]{Ours1} in the case of the figure-eight (with
   a different set of coordinates).  In the general case presented in \cite{Ours2}, the presence of
   powers of all the boundary lengths in \eqref{eq:form_coeff_phi} is a technical challenge, solved
   by expressing various boundary lengths as functions of our coordinate systems and studying their
   variations.
 \end{rem}

 \subsection{Generalized eights}
 \label{s:gene_eight}

 The proof that we sketched for the figure-eight generalizes to a whole family of topologies, which
 we call \emph{generalized eights}:
 \begin{defa} \label{d:GE} Let ${\curve}$ be a multi-loop filling a surface ${{\Sf}}$. Let $\cN$ be
   a regular neighbourhood of~${\curve}$ in ${{\Sf}}$. We say that ${\curve}$ is a \emph{generalized
     eight} if no connected components of $\partial \cN$ is homotopically trivial in ${{\Sf}}$
  \end{defa}

  \begin{figure}[h]
    \centering
    \includegraphics[scale=0.7]{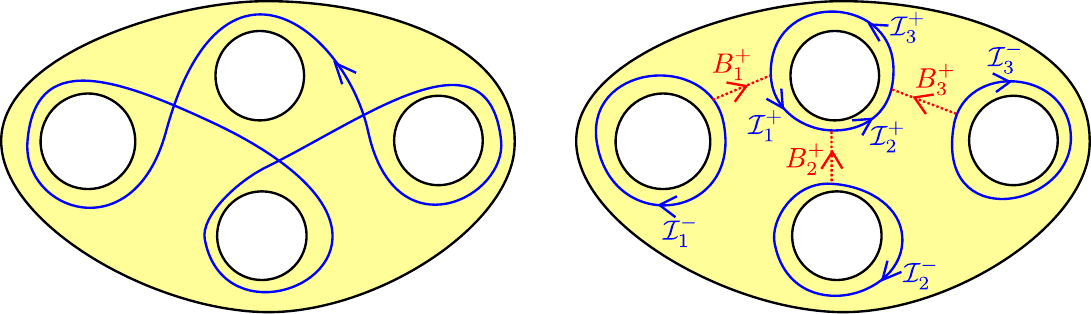}
    \caption{A generalized eight with $r=3$ intersections, and the segments $B_k^+$ and
      $\cI_j^\pm$ in its diagram.}
    \label{fig:generalized_eight}
  \end{figure}
  
  If a generalized eight ${\curve}$ fills a surface ${{\Sf}}$ of Euler characteristic $-r$, then
  ${\curve}$ has exactly $r$ intersection points. When we ``open'' these intersections, we obtain a
  \emph{diagram}, that is a multicurve $(\beta_1, \ldots, \beta_N)$ together with a family of $r$
  segments $B_1^+, \ldots, B_r^+$ corresponding to the intersections we removed -- see Figure
  \ref{fig:generalized_eight}.  Each of the curves $\beta_j$ and segments $B_k^+$ are naturally
  oriented by the original orientation of ${\curve}$.

 Conversely, a generalized eight ${\curve}$ is fully described, up to homotopy, by
 \begin{itemize}
\item an oriented multicurve $\beta=(\beta_1, \ldots, \beta_N)$ with $N$ components;
\item $r$ oriented segments $B_1^+,\ldots, B_r^+$, mutually disjoint and touching $\beta$ only at
  their endpoints, always leaving $\beta$ from its left and arriving at another point in $\beta$ on
  its right,
\item such that the union $\beta\cup B_1^+\cup,\ldots \cup B_r^+$ fills ${{\Sf}}$, with no disks in
  the complement.
\end{itemize}

From these data, a multi-loop ${\curve}$ can be recovered as follows. For $j \in \{1, \ldots, r\}$,
let $B_j^-$be the segment $B_j^+$ with opposite orientation.  For $\eps \in \{\pm\}$, call
$b(j, \eps)=t(B_j^{\eps})$ the terminus of the segment $B_j^{\eps}$, and $\cI_j^\eps$ the subsegment
of $\beta$ leaving from $b(j, \eps)$ and ending at the next point $b(j', \eps')$ along $\beta$.  Let
$p(j, \eps)=B_j^{\eps}\smallbullet \cI_j^\eps$ be the concatenation of the path $B_j^{\eps}$
followed by $\cI_j^\eps$. There is a unique way to concatenate together all the paths $p(j, \eps)$,
which gives the homotopy class of a generalized eight ${\curve}$.

Conversely, all generalized eights are obtained this way, up to homotopy.  We see that the segments
$B_1^+,\ldots, B_r^+$ capture the self-intersections of ${\curve}$. This description of the homotopy
class of ${\curve}$ is referred to as the \emph{diagram representing ${\curve}$}.

 \subsection{Coordinate system adapted to a generalized eight }
 \label{s:nc}
 
 In this section, we show how the bars $B_j^\pm$ and the segments $\cI_j^\pm$ can be used to define
 coordinates on the Teichm\"uller space~$\cT^*_{\g, \n}$, which are well adapted to expressing the
 geodesic length function~$Y \mapsto \ell_Y({\curve})$. This is essential to our analysis because an
 integration on the space~$\cT_{\g,\n}^*$, and the length function $\ell_Y(\curve)$, appear in the
 integral expression for $\av[\type]{F}$, Theorem \ref{thm:density_expr}.

 Note that the filled surface ${{\Sf}}$ has Euler characteristic $-r$, and thus $-r=2-2 \g-\n$.  The
 Teichm\"uller space $\cT^*_{\g, \n}$ then is a real-analytic manifold of dimension $3r$.

 Consider a point $(\x, Y)\in \cT^*_{\g, \n}$. We pick the representative in the Teichm\"uller space
 so that the components of $\beta$ are simple closed geodesics on~$Y$. We denote by
 $\overline{B}_j^\pm$ the orthogonal geodesic segments from $\beta$ to itself in the same homotopy
 class as the bar $B^\pm_j$ with endpoints gliding along $\beta$. As before, $\bar b(j, \eps)$
 denote the terminus of $\overline{B}_j^\eps$.  While the endpoint $b(j, \eps)$ glide to
 $\bar b(j, \eps)$, the segments $\cI_j^\eps$ glide to new subsegments $\overline{\cI}_j^\eps$ of
 $\beta$.

We then define the following coordinates:
\begin{itemize}
\item {\bf{Coordinates $L$.}} For $j \in \{1, \ldots, r\}$, we let $L_j=L_j(\x, Y)$ the (positive)
  length of $\overline{B}_j$.  
\item {\bf{Coordinates $\theta$.}} For $j \in \{1, \ldots, r\}$, $\eps=\pm$, we denote by
  $\theta_j^\eps=\theta_j^\eps(\x, Y)$ the {\em algebraic} length of $\overline{\cI}_j^\eps$:
  positive if it goes in the same direction as $\beta $, negative otherwise. 
\end{itemize}
We write $\vec{L}=(L_1, \ldots, L_r)$ and
$\vec{\theta}=(\theta_j^\eps)_{(j, \eps)\in \{1, \ldots, r\} \times \{\pm\}}$.  We can already note
that $(\vec{L}, \vec{\theta})$ belongs to a subset of $\R^{3r}$, a space of same
dimension as $\cT^*_{\g, \n}$.

For reasons that still remain to be explored, the Weil--Petersson measure on $\cT^*_{\g, \n}$ admits
a very nice expression, generalizing the one obtained in \eqref{e:magiccv}  (note that, for the
 figure-eight, $r=1$, $x_1, x_2$ are the $\theta$-coordinates and $T$ is the $L$-coordinate).
 \begin{prp}[{\cite{Ours2}}]\label{p:maindet}
   The map $(\x, Y)\in \cT^*_{\g, \n}\mapsto (\vec{L}, \vec{\theta})\in \R^{3r}$ is injective,
   and
   \begin{align}\label{e:maindet}\prod_{k=1}^{\n} \sinh\div{x_k}\d \x \d
     \mathrm{Vol}^{\mathrm{WP}}_{\g, \n, \x}(Y) = {2^r} \prod_{i=1}^{N}\sinh^2\div{y_i}
      \prod_{j=1}^r \sinh(L_j) \d L_j \d \theta_j^+ \d \theta_j^-
\end{align}
where, for any $1 \leq i \leq N$, $y_i$ denotes the length of the geodesic component $\beta_j$ of $\beta$ on $Y$.
 \end{prp}
 We notice, in particular, the appearance of the Friedman--Ramanujan functions $\sinh^2\div{y_j} $.
 
 The other nice feature of the $(\vec{L},\vec{\theta})$ coordinate system is that the length of
 ${\curve}$ admits a nice explicit expression which generalizes
 \eqref{e:eightformula}. 
 \begin{prp}
   There exists a family of vectors $\alpha^\delta\in \{-1,0, 1\}^{r}$, indexed by
   $\delta \in \{ \pm 1 \}^{2r}$ and depending on the diagram representing ${\curve}$, such
   that
   $$\ell_Y( {\curve})
   = 2\argcosh \Bigg( \sum_{\substack{\delta \in \{ \pm 1 \}^{2r} \\ \delta_1 \ldots \delta_{2r} =
         +1}} \cosh(\alpha^\delta \cdot \vec{L}) \hyp_{\delta} \div {\vec{\theta}} \Bigg)$$ where
   the notation $ \hyp_{\delta}$ is defined in \eqref{e:hypnotation}.

   Furthermore, $\alpha^{(1, \ldots, 1)}=(0, \ldots, 0)$, and there exists
   $\delta \in \{ \pm 1 \}^{2r}$ such that $\alpha^{\delta}=(1, \ldots, 1)$.
  \end{prp}
  Because there exists $\delta \in \{ \pm 1 \}^{2r}$ with
  $\cosh(\alpha^\delta \cdot \vec{L}) =\cosh (L_1+\ldots+L_r)$, we can consider that the $L_i$
  contribute twice to $\ell_Y( {\curve})$. This will allow to treat $ \vec{L}$ as ``neutral
  parameters'' in our analysis, meaning that we will not need to apply the operator $\cL$ to these
  variables in order to cancel their exponential behaviour (as was the case with the variable $T$ in
  the case of the figure-eight, c.f.  Remark~\ref{rem:neutral_eight}).

 Using Theorem \ref{thm:density_expr}, we are able to generalize Lemma
  \ref{e:Vconvol} and express the volume function as a linear combination of functions of the form
  \begin{align*}
   \int  f_1\star  f_2\star\ldots \star f_{2r}|^{h_{\vec{L}} }_{\varphi_{\vec{L} }}\,
    \prod_{j=1}^r \sinh(L_j)\d L_j, 
  \end{align*}
  where:
  \begin{itemize}
  \item the $(f_j)_{j}$ are Friedman--Ramanujan functions (some of them might actually be Dirac
    delta distributions, as in \eqref{e:handle}, in which case the integration is simply reduced to
    a smaller domain);
\item for a fixed $ \vec{L}$, the
  function
\begin{equation}
h_{\vec{L}} : \vec{\theta}\mapsto 2\argcosh \left(  \sum_{\substack{\delta \in \{ \pm 1 \}^{2r} \\ \delta_1 \ldots \delta_{2r} = +1}} \cosh(\alpha^\delta \cdot  \vec{L})   \hyp_{\delta} \div {\vec{\theta}} \right)
\end{equation}
appears to fall exactly into the class defined in Notation~\ref{nota:hL};
\item the function $\varphi_{\vec{L}}$ contains several defects: the indicator function of the domain
  of the coordinates $(\vec{L},\vec{\theta})$ as well as the boundary lengths mentioned in
  Remark~\ref{rem:boundary_comp}.
\end{itemize}
This allows a similar treatment as the one we sketched for the figure-eight. A caveat is that,
unless additional hypotheses are made on the generalized eight (the ``non-crossing hypothesis'' in
\cite{Ours2}), some of the variables $\theta_j^\eps$ vary in an interval allowing negative values.
If this is the case, Proposition \ref{p:mainclass} does not apply directly to the function
$h_{\vec{L}} $; this requires additional technical work, that we do not describe here.

\subsection{One word about other topologies}

Generalized eights form a family of loops that is simple to describe and enumerate. Beyond
these cases, it may seem daunting to study every possible topology. Actually, generalized
eights are the worse possible scenario for the proof of Theorem \ref{t:leviathan}, and any other situation is ``easier'', as we shall
here briefly explain.

We have hinted in Remark \ref{rem:neutral_eight} and \S \ref{s:nc} that there are two kinds of
variables in our new set of coordinates on the Teichm\"uller space: \emph{active} parameters, for
which we need to apply the operator $\cL$ to cancel exponential terms, versus \emph{neutral}
parameters, which do not cause any exponential behaviour and can be controlled by simple
inequalities.  For generalized eights, these are respectively the parameters $\vec{\theta}$ and
$\vec{L}$.  Actually, for more complicated topologies, we can define similar parameters
$\vec{\theta}$ and $\vec{L}$, but some of the $\vec{\theta}$ coordinates will become ``neutral parameters''.

Indeed, in general, if a loop $\curve$ fills a surface ${\Sf}$ endowed with a metric
$(\x, Y)\in \cT^*_{(\vg, \vn)}$, then we always have
\begin{equation}
  \label{eq:comparison_basic}
  \sum_{i=1}^{n_\Sf} x_i \leq 2\ell_Y(\curve).
\end{equation}
This trivial inequality is obtained by observing that the length of the boundary of a regular
neighbourhood of $\curve$ can be made arbitrarily close to $2\ell_Y(\curve)$. Interestingly, if we
were to manage to remove the factor $2$ and prove the improved inequality
\begin{equation}
  \label{eq:comparison_double_fill}
  \sum_{i=1}^{n_\Sf} x_i \leq \ell_Y(\curve)
\end{equation}
then it would be straightforward to deduce, using simple upper bounds, that the volume function
$V_g^{\eqc{\Sf,\curve}}$ is a Friedman--Ramanujan remainder, i.e. an element of $\cR_w$. In
other words, the goal of this section, Theorem \ref{t:leviathan}, becomes trivial. This is the
reason why the contribution of non-primitive geodesics in unimportant in our trace method, as proven
in Lemma \ref{lem:bound_r1}. See \cite[Proposition 8.6]{Ours1} in the case of loops filling a pair
of pants.

In \cite{Ours1,Ours2}, we introduce the following concepts.
\begin{defa}
  Let $\curve$ be a loop filling a surface $\Sf$. A simple portion $\cI$ of $\curve$ is an open
  subsegment of $\curve$ that does not contain any self-intersection. The simple portion $\cI$ is said
  to be \emph{shielded} if it runs along a contractible component of $\Sf \setminus \curve$. We say
  that $\curve$ is \emph{double-filling} if all of its simple portions are shielded.
\end{defa}
We then refine the inequality \eqref{eq:comparison_basic} and rather write
\begin{equation*}
  \sum_{i=1}^{n_\Sf}x_i \leq \ell_Y(\curve) + \sum_{\cI \text{ non-shielded}} \ell(\cI)
\end{equation*}
where the sum runs over all maximal simple portions of $\curve$ that are non-shielded.  This formula
become \eqref{eq:comparison_basic} in the case of generalized eights, because all portions are
non-shielded. However, in the case of double-filling loops, it reads
\eqref{eq:comparison_double_fill}, meaning Theorem \ref{t:leviathan} is trivial (and hence every
single parameter is neutral). In \cite{Ours2}, we explain how to deal with intermediate cases,
and apply the operator $\cL$ to all active parameters.

   \section{Tangles and the Moebius inversion formula}\label{s:tangles}
   \label{s:tangles_moebius}
   
   In this section, we shall present our solution to the issues raised in \S \ref{s:tangle_pb},
   caused by the presence of tangles. These results are presented in the standalone article
   \cite{Moebius-paper}.

   We first define the notion of tangle in \S \ref{s:deftangles}. Then, Theorem \ref{cor:TF_curves} states
   that in $(\kappa, \tfL)$-tangle free surfaces of a given Euler characteristic $-\chi$, with $\tfL$
   of magnitude $\log g$, the number of possible topologies for a periodic geodesic of length
   $L=A\log g$ is at most polynomial in $L$ (instead of exponential). This means that removing our
   tangles does remedy the issue of exponential proliferation of topologies presented in \S
   \ref{s:tangle_pb}.

   We then present in \S \ref{s:moebius} a Moebius inversion formula, inspired by Friedman's work
   \cite{friedman2003}, which is a form of inclusion-exclusion argument allowing to adapt our trace
   method to the new probability model conditioned on the set of tangle-free surfaces.

   \subsection{c-surfaces and extended notions of moduli spaces}

   Before we define the notion of tangles, let us set a few notations which will be useful in the following.
   
   We order the pairs $(g, n)\in \N_0^2$ by the order relation: $(g, n)\leq (g', n')$ if
   $2g'+n'-2 > 2g+n-2$, or if $2g'+n'-2 = 2g+n-2$ and $g'\geq g$.
   For $k\geq 1$, denote by ${\mathfrak{S}}^{(k)}$ the set of $2k$-tuples
   $(\vg, \vn)=((g_1, n_1), \ldots, (g_k, n_k))\in (\N_0^2)^k$, satisfying $2-2g_i-n_i\leq 0$ and
   $(g_i, n_i)\leq (g_{i+1}, n_{i+1})$.
 
   \begin{defa}[{\cite[Definition 2.1]{Moebius-paper}}] A c-surface of signature
     $(\vg, \vn)\in {\mathfrak{S}}^{(k)}$ is a non-empty topological space $\Sf$ with $k$ connected components
     $\Sf= \sqcup_{j=1}^k \tau_j$ labelled from $1$ to $k$, where:
     \begin{itemize}
     \item $\tau_i$ is a $1$-dimensional oriented manifold diffeomorphic to a circle if
       $(g_i, n_i)=(0, 2)$;
     \item  if $2-2g_i-n_i<0$, $\tau_i$ is a $2$-dimensional orientable manifold with (or without) boundary, of signature $(g_i, n_i)$, with a numbering of boundary components by the integers $1, \ldots, n_i$.
     \end{itemize}
   \end{defa}

   The c-surface $\Sf$ can be decomposed as $\Sf=(c, \sigma)$, where $c =(c_1, \ldots, c_j) $ and
   $\sigma =(\sigma_1, \ldots, \sigma_m)$ are, respectively, collections of $1$-dimensional
   manifolds and $2$-dimensional orientable compact manifolds with boundary ($k=j+m$ is the total
   number of connected components). 
   We shall denote by $|\Sf|=j$ the number of 1d components, and $\chi(\Sf)$ the total absolute Euler
   characteristic, that is $\chi(\Sf)=\sum_{i=1}^m \chi(\sigma_i)$ where
   $\chi(\sigma_i)=2g_{i+j}+n_{i+j}-2$.

   The motivation for including the case $\chi = 0$ in the definition of c-surface is that the
   tangles we wish to exclude correspond to several geometric patterns: certain embedded
   subsurfaces, but also short closed geodesics.

   We extend the notations $\cT_{g,n}^*$ and $\cM_{g,n}^*$ to the cases $g=0$ and $n=2$,
   corresponding to $2-2g-n=0$: $\mathcal{T}_{0,2}^*= \mathcal{M}_{0,2}^*$ is the space of
   Riemannian metrics on an oriented circle, modulo orientation preserving isometry.  If
   $Z\in \mathcal{M}_{0,2}^*$, we denote by $\ell(Z)$ the length of $Z$. The map $Z\mapsto \ell(Z)$
   allows to identify $\mathcal{M}_{0,2}^*$ with $\R_{>0}$.

   For $(\vg, \vn)=((g_1, n_1), \ldots, (g_k, n_k))\in {\mathfrak{S}}^{(k)}$, we define
   $$\mathcal{M}_{(\vg, \vn)}^* := \prod_{i=1}^k \mathcal{M}_{g_i,n_i}^*.$$ This can be interpreted as
   the moduli space of Riemannian metrics on a c-surface with $k$ labelled connected components, of
   respective signatures $(g_1, n_1), \ldots, (g_k, n_k)$.  Finally define
 $$\mathbf{M} := \bigsqcup_{k=1}^{+\infty}\bigsqcup_{(\vg, \vn)\in  {\mathfrak{S}}^{(k)}}\mathcal{M}_{(\vg, \vn)}^*.$$
 The space $\mathbf{M}$ is the space of Riemannian metrics on all possible c-surfaces, hyperbolic on the 2d-components, modulo isometries preserving the orientation and the numbering of connected and boundary components.
 
 \subsection{Definition of tangles} \label{s:deftangles}
  
 In order to overcome the difficulties exposed in \cref{s:tangle_pb}, we need to remove some ``bad''
 surfaces containing many geodesics of length $A\log g$. Let us introduce a notion of
 \emph{tangles}, associated with two parameters $\kappa, \tfL >0$ with $\kappa< \tfL$.

\begin{defa}[{\cite[Definition 2.4]{Moebius-paper}}]
  \label{def:tangle}
  We call $Z\in \mathbf{M}$ a {\em{tangle}} (or $(\kappa, \tfL)$-tangle if we wish to keep explicit the dependency on $(\kappa, \tfL)$) if
   \begin{itemize}
  \item either $Z\in \mathcal{M}_{0,2}^*$ and $\ell(Z)\leq \kappa$;
  \item or $Z\in \mathcal{M}_{g,n}^*$ with $2-2g-n=-1$, and $\ell^{\mathrm{max}}(\partial Z) \leq \tfL.$
   \end{itemize}
  \end{defa}

  As a consequence, tangles are either short simple closed geodesics (of length $\leq \kappa$) or
  embedded pairs of pants or once-holed tori with short geodesic boundary (of total length $\leq
  \tfL$). We know that both those scenarios are linked with difficulties in the study of spectral gaps:
  Lipnowski--Wright performed an inclusion-exclusion in \cite{lipnowski2021} to remove the former,
  whilst we exposed in \S \ref{s:tangle_pb} the issues related to the latter.

\begin{nota}
  We denote as $\Atf$ the set of hyperbolic surfaces of genus $g$ that do not contain
  $(\kappa, \tfL)$-tangles, and call them \emph{tangle-free surfaces}.
\end{nota}

 For $\tfL=\kappa \log g$, the probability of $\Atf$ can be estimated as follows:

\begin{lem}[{\cite[Theorem 4.2]{mirzakhani2013} and \cite[Theorem 5]{monk2021a}}]
  \label{lem:prob_TF}
  For any small enough $\kappa > 0$, any large enough $g$, $\tfL= \kappa\log g$, 
  \begin{equation}\label{e:prob_TF}
    1 - \Pwp{\Atf} = \O{\kappa^2 + g^{\frac 3 2 \kappa-1}}.
  \end{equation}
\end{lem}

We will always take $\kappa < \min(2/3, 2\argsh 1)$ -- in fact, we will later need to take $\kappa$
even smaller, and depending on $A$. In this case, the probability of $\Atf$ goes to one as
$g \rightarrow + \infty$, and it is therefore reasonable to condition on this subset to prove
\eqref{e:dream}.

Taking $\tfL=\kappa\log g$ and $L=A\log g$, we can write
\begin{equation*}
  \Pwp{\delta \leq \lambda_1 \leq \frac 1 4 - \alpha^2 - \epsilon} 
  \leq    \Pwp{\left\{\delta \leq \lambda_1 \leq \frac 1 4 - \alpha^2 - \epsilon\right\}\cap \Atf} + 
  1-  \Pwp{\Atf}.
\end{equation*}
This allows us to modify the statement of Lemma~\ref{lem:selberg_reformulated}, and replace it by:
  \begin{equation}\label{e:conditioned}
    \begin{split}
      &   \Pwp{\delta \leq \lambda_1 \leq \frac 1 4 - \alpha^2 - \epsilon} 
      \\ & \leq    \frac{C}{g^{(\alpha + \epsilon) A}} \Big(
           \avbtf{\ell \,e^{-\ell/2} \cD^m H_L(\ell)}
           + g \log(g)^2\Big)  + 
            C(\kappa^2 + g^{\frac 3 2 \kappa-1})
    \end{split}
    \end{equation}
    where $\chic = 3A$  and we define
   \begin{align}  
  \avbtf{F} 
  \label{e:reduce}  
     :=\Ewp{ \sum_{\type : \chi(\type) \leq \chic} \sum_{\substack{\gamma\sim \type}} F(\ell(\gamma))
     \, \bbbone_{\Atf}(X)}.
 \end{align}
 Note that we have used Proposition \ref{p:apriori} to reduce the sum to filling types of bounded
 Euler characteristic.

 \begin{summ}
   We have introduced a subset of tangle-free surfaces $\Atf$ by which we mean to condition the
   Weil--Petersson probability model, to remove the issues presented in \S \ref{s:tangle_pb}. It is
   straightforward to adapt the trace method and condition by $\Atf$ -- this leads to now
   considering the conditioned average $\avbtf{\ell \,e^{-\ell/2} \cD^m H_L(\ell)}$.
 \end{summ}

\subsection{Counting topological types of geodesics in tangle-free surfaces}
 
Let us show how removing tangles solves the problem of exponential proliferation of local types
presented in \S \ref{s:tangle_pb}.

\begin{defa} For $0<\kappa< \tfL$, $L>0$ and $\Sf$ a filling type, we denote by
  $\mathrm{Loc}_{ \Sf}^{\kappa,\tfL, L}$ the set of local topological types $\eqc{\Sf, \curve}$, such
  that there exists a $(\kappa, \tfL)$-tangle-free surface $Z$ and a homeomorphism $\phi: \Sf\To Z$
  with $\ell_Z(\phi(\curve))\leq L$.  If $\chic$ a positive integer, we let
$$\mathrm{Loc}_{\chic}^{\kappa,\tfL, L} := \bigcup_{\Sf, \chi(\Sf)\leq \chic}\mathrm{Loc}_{\Sf}^{\kappa,\tfL, L}$$ 
where the union is over the (finite) set of filling types $\Sf$ with $\chi(\Sf)\leq \chic$.
\end{defa}

By definition, we can write for any test function $F$
\begin{equation}
  \label{eq:TF_Loc}
  \avbtf{F}
  = \Ewp{ \sum_{\type : \chi(\type) \leq \chic} \sum_{\substack{\gamma\sim \type}} F(\ell(\gamma))
     \, \bbbone_{\Atf}(X)}
  = \sum_{\type \in \mathrm{Loc}_{\chic}^{\kappa,\tfL, L}} \avTbtf{F}
\end{equation}
where, for any $\type \in \mathrm{Loc}_{\chic}^{\kappa,\tfL, L}$,
\begin{align} \label{e:avtbtf}  \avTbtf{F} :=\Ewp{\sum_{\substack{\gamma\sim \type}} F(\ell(\gamma)) \, \bbbone_{ \Atf}(X)}.
 \end{align}
We prove the following upper bound on the cardinality of $\mathrm{Loc}_{\chic}^{\kappa,\tfL, L}$ at
logarithmic scales in $g$.

\begin{thm}[{\cite[Theorem 3.5]{Moebius-paper}}]
  \label{cor:TF_curves}
  Let $0 < \kappa < 1$, $\tfL=\kappa \log g$, $A \geq 1$, $L = A \log g$ and $\chic \geq 1$. Then,
    \begin{equation}
    \label{eq:bound_number_type_TF}
    \# \mathrm{Loc}_{\chic}^{\kappa,\tfL, L} = \O[\kappa, A,\chic]{(\log g)^{\beta_{\kappa,A, \chic}}}
  \end{equation}
  for a $\beta_{\kappa, A, \chic} > 0$ depending only on $\kappa$, $A$ and $\chic$.\end{thm}

\begin{summ}
  At the scale $\tfL=\kappa \log g$ and $L=A\log g$, conditioning the Weil--Petersson measure by the
  set of tangle-free surfaces does circumvent the exponential proliferation exhibited in \S
  \ref{s:tangle_pb}. Noticeably, the last sum in \eqref{eq:TF_Loc} is reduced to a polynomial number
  of terms, rather than an exponential one. It follows that we can aim to bound each term of this
  sum individually, and obtain a decent bound on the overall quantity
  $\avbtf{\ell \,e^{-\ell/2} \cD^m H_L(\ell)}$.
\end{summ}

However, the Friedman--Ramanujan argument from \S \ref{s:topo} is affected by the apparition of an
indicator function $\1{\Atf}$ in the average. The aim of the rest of this article is to explain how
to adapt the integration-by-parts argument to the conditioned average, and hence conclude to the
proof of \eqref{e:dream}.

 \subsection{Moebius inversion formula}
 \label{s:moebius}

 We construct a ``generalized Moebius inversion formula'', allowing to replace the indicator function
 $\bbbone_{\Atf}$ by an alternating sum of counting functions. The terminology comes from J. Friedman
 \cite[Proposition 9.5]{friedman2003}, in a similar context for random regular graphs.

 \subsubsection{Sub-c-surfaces and their geodesic representatives}
 \begin{defa}\label{d:subsurface}
   If $\Sf$ is a c-surface, we call {\em{sub-c-surface}} of $\Sf$ a subset $\Sf'\subset \Sf$
   which (for the induced topology and differential structure)
   is a c-surface, and such that
   \begin{itemize}
   \item the 1d components of $\Sf'$ are homotopically non-trivial in $\Sf$;
   \item the boundary curves of the 2d components of $\Sf'$ are homotopically non-trivial in $\Sf$, and are pairwise non-homotopic.
   \end{itemize}
   In other words, we ask that the boundary $\partial \Sf'$ is a multicurve in $\Sf$.
 \end{defa}

 Let $\Sf$ be a c-surface of signature $(\vg, \vn)$ and $\Sf' = (c',\sigma')$ be a sub-c-surface of $\Sf$.
 If $\Sf$ is endowed with a metric $ Z\in \cT^*_{(\vg, \vn)}$, then $\Sf'$ is isotopic to a geodesic
 sub-c-surface $(c'_Z, \sigma'_Z)$ (that is to say, a sub-c-surface with geodesic boundary),
 called the {\em{geodesic representative}} of $\Sf'$ in $Z$. 

We denote by $\cS(Z)$ the set of geodesic sub-c-surfaces of $\Sf$ endowed with the metric $Z$. 

\subsubsection{Derived tangles and maximal tangle}

A challenge when performing inclusion-exclusion arguments to remove tangles is that we will need to
enumerate families of tangles, and hence detail how different tangles might combine.  This motivates
the following definition, directly inspired by the work of Friedman \cite{friedman2003}.
  
  \begin{defa}[{\cite[Definition 2.5]{Moebius-paper}}]
  \label{def:dtangle}
  We call $Z\in \mathbf{M}$ a {\em{derived tangle}} (or $(\kappa, \tfL)$-derived tangle) if it is
  filled by tangles. In other word, $Z\in \prod_{i=1}^k \mathcal{M}_{g_i,n_i}^*$ has $k$ connected
  components, some of them circles and some of them hyperbolic manifolds with geodesic boundary; and
  there is a sequence of submanifolds $ (Z_n)_{n\geq 1}$ which are either periodic geodesics, or
  surfaces of Euler characteristic $-1$, which are tangles in $Z$, and whose union fill $Z$.

  Denote by ${\mathbf{D}}_{\kappa, \tfL} \subset \mathbf{M}$ the set of all $(\kappa, \tfL)$-derived
  tangles.
\end{defa}
In the statement above, the notion of subset filling a c-surface is defined by naturally extending
the notion of subset filling a surface (see \cite[Definition 2.2]{Moebius-paper}).

It is clear that any hyperbolic surface $X$ contains a maximal derived tangle, defined as the
surface with geodesic boundary filled by all the tangles contained in $X$. We denote it by
$\tau_{\kappa,\tfL}(X)$.  Note that $\tau_{\kappa,\tfL}(X)$ may not be a \emph{sub-c-surface} of
$X$, because several of the boundary components of $\tau_{\kappa,\tfL}(X)$ may be the same in $X$.

 \subsubsection{The ``Moebius function''}
 \label{s:Moebius}
 Let $\kappa, \tfL$ be two real numbers. Assume as before (mostly for simplicity) that
 $\kappa<\tfL$, and that $\kappa <2\argsh(1)$, so that all closed geodesics of length $\leq \kappa$
 in a hyperbolic manifold are simple and pairwise disjoint \cite[Theorem 4.1.6]{buser1992}.
 
 \begin{thm}[{\cite[Theorem 4.3]{Moebius-paper}}] \label{t:moebius}There exists a (unique) function
   $\mu = \mu_{(\kappa, \tfL)} : {\mathbf{M}} \To \R$, such that:
 \begin{itemize}
  \item the restriction of $\mu$ to each $\cM^*_{(\vg, \vn)}$ is invariant under the action of the full diffeomorphism group ({\emph{i.e.}} possibly permuting boundary components and connected components);
\item for $Z\in {\mathbf{M}} $, $\mu(Z)=0$ if $Z\not \in {\mathbf{D}}_{\kappa, \tfL} $;
\item for every $Z\in {\mathbf{D}}_{\kappa, \tfL} $,
\begin{align}\label{e:fundidentity}
1= \sum_{\tau \in \cS(Z)} \mu(\tau),
\end{align}
where $\mu(\tau)$ is defined as $\mu(\bar{\tau})$, where $\bar{\tau}\in {\mathbf{M}}$ is the
equivalence class of $\tau$.
 \end{itemize}
 In addition, the function $\mu$ we construct satisfies the following:
 \begin{itemize}
 \item if $(\vg, \vn)=(\underbrace{(0,2), \ldots, (0,2)}_{j})$, and $Z=(c_1, \ldots, c_j)\in \cM^*_{(\vg, \vn)}$, then
   \begin{align}\label{e:circles}
     \mu(Z)=
     \begin{cases}
       \frac{(-1)^{j+1}}{2^j j!} & \text{if } \ell^{\mathrm{max}}(c)\leq \kappa,\\
       0 & \text{otherwise;}
     \end{cases}
   \end{align}
 \item if $Z=(c, \sigma)$ is the decomposition of $Z$ into 1d and 2d part, then the function $\mu$
   has the multiplicative property $\mu(Z)= -\mu (c)\mu(\sigma)$;
\item there exists explicit increasing sequences $U(n)$ and $V(n)$ such that
 \begin{align}\label{e:uppermu}|\mu(Z)|\leq U(\chi(Z)) \, e^{\tfL V(\chi(Z))}.
 \end{align} 
 \end{itemize}
\end{thm}

We then obtain, as a straightforward consequence, the following formula, called the Moebius
inversion formula.

\begin{cor}[Moebius inversion formula {\cite[Corollary
    1.3]{Moebius-paper}}] \label{p:moebius-goal}For every $g$, for every compact hyperbolic surface
  $X$ of genus $g$,
    \begin{align} \label{e:inversionformula}1- \bbbone_{\Atf}(X)= \sum_{\tau\in
        \cS(X)} \mu(\tau).
\end{align}
 \end{cor}

 \section{Obtention of the spectral gap}
 \label{s:obtention}

 Let us now explain how to use the Moebius inversion formula to show that a variant of Theorem
 \ref{thm:exist_asympt_type_intro} applies to the conditioned averages in
 \eqref{e:conditioned},  and conclude to the proof of \eqref{e:dream}.

 \subsection{Application of the Moebius inversion formula}
\label{sec:appl-moeb-invers}

Before applying the Moebius formula, it will be convenient to introduce the following event.
\begin{nota}
  For an integer $Q$ which shall be specified later, we denote as ${\MC}_X(Q)$ the set of
  multicurves which separate $X$ into at most $Q$ connected components, and let
  \begin{align*}\cB_g^{ \kappa, Q}
    :=
    \left\{ X\in \cM_g \, : \,
    \forall c  \text{ multicurve s.t. } \ell^{\mathrm{max}}(c) \leq \kappa, 
    \, c \in \MC_X(Q)
    \right\}.
    \end{align*}
\end{nota}
By definition, we have that $\Atf \subset\mathcal{B}_g^{\kappa, Q}$ for any $Q$, and hence the
probability of $\mathcal{B}_g^{\kappa, Q}$ goes to~$1$ as $g \rightarrow + \infty$, and
$\1{\Atf} = \1{\Atf} \1{\mathcal{B}_g^{\kappa, Q}}$. Let us further condition the Weil--Petersson
model by the set $\cB_g^{ \kappa, Q}$, which will allow us to limit the possibilities when
enumerating the 1d-part of tangles.  The Moebius inversion formula allows to rewrite
\eqref{e:avtbtf} the following way.
 
 \begin{prp}[{\cite{Ours2}}]
   \label{prp:proba_small_Moebius}
  For $L = A \log g$, $\tfL = \kappa \log g$, $\chitau = 2A+4$, $Q = 2A + 5+6\chitau$ and provided
  $2 \kappa   V(\chitau) + 12 \chitau \kappa < 1$, with $V$ the sequence in \eqref{e:uppermu}, we have:
   \begin{align} \label{e:muav}
    \avTbtf{F} = 
    \avb[\type]{F} 
     -\Ewpo \Bigg[ \sum_{\substack{\tau = (c, \sigma)\in \cS(X) \\ \chi(\tau) <\chitau \\ c\in \,
   {\MC}_X(Q)}}
   \sum_{\substack{\gamma \in \geod(X)\\ \gamma\sim \type}} 
   \mu_{\kappa,\tfL}(\tau)  F(\ell(\gamma)) \Bigg] + \O[A]{1}.
 \end{align}
\end{prp}
The proof of Proposition \ref{prp:proba_small_Moebius} consists in bounding the remainder in
\eqref{e:muav} above by
\begin{align}
   \Ewp{(1-\1{\mathcal{B}_g^{\kappa,Q}}(X) + \1{\chi(\tau_{\kappa,\tfL}(X)) \geq \chitau})
    \Bigg(1 +  \sum_{\substack{\tau = (c, \sigma)\in \cS(X) \\ \chi(\tau) <\chitau \\ c\in \,
   {\MC}_X(Q)}} |\mu_{\kappa,\tfL}(\tau)|\Bigg) \sum_{\substack{\gamma \in \geod(X)\\ \gamma\sim
  \type}}  |F(\ell(\gamma))| },
\end{align}
where we recall that $\tau_{\kappa,\tfL}(X)$ denotes the maximal $(\kappa,\tfL)$-derived tangle of
$X$.  We conclude by using the bound \eqref{e:uppermu} on the Moebius function $\mu_{\kappa,\tfL}$
and a Cauchy--Schwarz argument, relying on the following probabilistic bounds on the events
$\mathcal{B}_g^{\kappa, Q}$ and $\chi(\tau_{\kappa,\tfL}(X)) \geq \chitau$:
\begin{align}
  &1-\Pwp{ \cB_g^{ \kappa, Q}} =\cO_{\kappa, Q} \Big( \frac{1}{g^{Q-1}}\Big) \\
  &\Pwp{\chi(\tau_{\kappa,\tfL}(X)) \geq \chitau} = \O[\chitau]{\frac{(\chitau \omega)^{c_\chitau}e^{3
    \chitau \omega}}{g^\chitau}}.
\end{align}

Looking back at \eqref{e:conditioned}, we recall that we aim at applying \eqref{e:muav} to
$F(\ell)= \ell e^{-\ell/2} \cD^m H_L(\ell)$.  We see that the first term $\avb[\type]{F}$ is already
dealt with in Proposition \ref{cor:FR_implies_small}.  Our main concern now is the term involving
the Moebius function $\mu_{\kappa,\tfL}$~: we shall show that a variant of Proposition~\ref{cor:FR_implies_small}
applies to it.

\subsection{Local topological types of lc-surfaces.}
\label{s:lcsurface}

Let us enrich the notion of local topological type, to a notion designed to sort the new terms
appearing in \eqref{e:muav}. The topological object to consider is now the pair formed by
$(\curve, \tau)$, where $\curve$ is a loop and $\tau$ is a c-surface. Such a pair will be called a
lc-surface. We define topological types of such objects in \cite[\S 4.2]{Ours1}.
  
   \begin{nota}\label{n:FTlc} To any $k\geq 0$ and any $(\vg, \vn)\in {\mathfrak{S}}^{(k)}$
   we shall associate a \emph{fixed} smooth
  oriented c-surface $\Sf$ of signature $(\vg, \vn)$. The data of $(\vg, \vn)$, or equivalently of the c-surface $\Sf$, is called a \emph{c-filling
    type}.
\end{nota}

\begin{defa}
  \label{def:local_type}
  A \emph{local lc-surface} is a triplet $(\Sf,\curve, \tau)$, where $\Sf$ is a c-filling type,
  $\curve$ is a loop in $\Sf$, $\tau$ is a sub-c-surface of $\Sf$, and the pair $(\curve, \tau)$
  fills $\Sf$. Two local lc-surfaces $(\Sf,\curve, \tau)$ and $(\Sf',\curve', \tau')$ are said to be
  \emph{locally equivalent} if $\Sf=\Sf'$ (i.e.  $\mathbf{g}_\Sf = \mathbf{g}_{\Sf'}$ and
  $\mathbf{n}_\Sf=\mathbf{n}_{\Sf'}$), and there exists a positive homeomorphism
  $\psi : \Sf \rightarrow \Sf$, possibly permuting the boundary components and  connected
  components of $\Sf$, such that $\psi \circ \curve$ is freely homotopic to $\curve'$ and
  $\psi(\tau)=\tau'$.  This defines an equivalence relation $\eq$ on local lc-surfaces. Equivalence
  classes for this relation are denoted as $\eqc{\Sf, \curve, \tau}$ and called \emph{local
    topological types} of lc-surfaces.
\end{defa}
 
  \begin{defa}
    Let $\eqc{\Sf, \curve, \tau}$ be a local topological type.

    For a pair $(\gamma, \tau_0)$ formed by a loop $\gamma$ and a sub-c-surface $\tau_0$ on a
    hyperbolic surface~$X$, we denote as $\Sf(\gamma, \tau_0)$ the surface filled by
    $(\gamma, \tau_0)$ in $X$ (defined in details in \cite[Definition 4.1]{Ours1}). Then,
    $(\gamma, \tau_0)$ is said to \emph{belong to the local topological
      type}~$\eqc{\Sf, \curve, \tau}$ if there exists a positive homeomorphism
    $\phi : \Sf(\gamma, \tau_0) \rightarrow \Sf$ such that the loops $\phi ( \gamma)$ and $\curve$
    are freely homotopic in~$\Sf$, and such that $\phi(\tau_0)=\tau$. In that case, we write
    $(\gamma, \tau_0) \sim \eqc{\Sf, \curve, \tau}$.
\end{defa}

\begin{defa} \label{e:def_Tan} Denote by $\LocTFlc$ the set of local topological types
  $\eqc{\Sf, \curve, \tau}$ of lc-surfaces of filling type $\Sf$, such that:
  \begin{itemize}
  \item $\eqc{\Sf(\curve),\curve} \in \mathrm{Loc}_{\chic}^{\kappa,\tfL, L}$, where
    $\Sf(\curve) \subset \Sf$ is the surface filled by $\curve$;
  \item $\chi(\tau)<\chitau$;
  \item there exists a hyperbolic surface $Z$ and a homeomorphism $\phi : \Sf\To Z$ with
    $\ell_Z(\phi(\curve))\leq L$ and such that the geodesic representative of $\phi(\tau)$ is a
    $(\kappa, \tfL)$-derived tangle in $Z$.
  \end{itemize}
\end{defa}

Let us define $\cS_Q(X)\subset \cS(X)$ as the set of geodesic subsurfaces of $X$, with a 1d-part
belonging to ${\MC}_X(Q)$.  In light of \eqref{e:muav}, the quantity we wish to study may be
rewritten as
\begin{align}\label{e:tensor}
  \sum_\Sf\sum_{{\type}\in \LocTFlc}
  \Ewpo \Bigg[ \sum_{\substack{ \gamma \in \mathcal{P}(X) \\ \tau\in \cS_Q(X)\\ (\gamma, \tau)\sim {\type}}} 
  \mu_{\kappa, \tfL} (\tau) F( \ell(\gamma) ) \Bigg]
  = \sum_\Sf \sum_{{\type}\in \LocTFlc }
  \avN[{\type}]{\mu_{\kappa, \tfL} \otimes F},
 \end{align}
 where now we let for any bounded function $\nu$ on $\mathbf{M}$
 \begin{align}  \label{e:total}
   \avN[\type]{\nu \otimes F}
   :=   \Ewpo \Bigg[ \sum_{\substack{ \gamma \in \mathcal{P}(X) \\\tau\in \cS_Q(X)\\ (\gamma, \tau)\sim {\type}}} 
 \nu(\tau) F( \ell(\gamma) ) \Bigg].
 \end{align}
 
 The first sum in \eqref{e:tensor} runs over the infinite set of all filling types $\Sf$. However,
 arguing like for~\eqref{e:bound_chi}, if $L=A\log g$ and $\tfL=\kappa\log g$, we may find
 $\chictau=\chictau(A, \chitau, Q)$ such that, for all $\kappa< 2/3$,
 \begin{align}\label{e:chi_reduc} \sum_\Sf \sum_{{\type}\in \LocTFlc }   \avN[{\type}]{\mu_{\kappa, \tfL} \otimes F}
   =\sum_{\substack{\Sf \\ \chi(\Sf)< \chictau}} \sum_{{\type}\in \LocTFlc }
   \avN[{\type}]{\mu_{\kappa, \tfL} \otimes F} +\cO_{A, \chitau, Q}(1).
 \end{align}

 One should remark that even if we fix the Euler characteristic of the c-surface $\Sf$, the sum
 $\sum_{\Sf, \chi(\Sf)< \chictau}$ still has infinitely many terms, because the c-surface $\Sf$ may have
 arbitrarily large number of 1d-components $|\Sf|$.  However, there is an obvious bijection from the
 set of c-surfaces with $|\Sf|=0$ to the set of c-surfaces with $|\Sf|=j$, which consists in adding $j$
 1-dimensional connected components.  Thus we can write
 \begin{align} \label{e:totalj}
 \sum_{\substack{\Sf \\ \chi(\Sf)< \chictau}} \sum_{{\type}\in \LocTFlc }   \avN[{\type}]{\mu_{\kappa, \tfL} \otimes F}
   =\sum_{\substack{\Sf \\ \chi(\Sf)< \chictau \\ |\Sf|=0}} \sum_{{\type}\in \LocTFlc}
   \sum_{j=0}^{+\infty}  \avN[{R_j\type}]{\mu_{\kappa, \tfL} \otimes F},
\end{align} 
where the operation $R_j$ consists in adding $j$ 1-dimensional connected components to a local
topological type.
 
\subsection{Adaptation of our arguments  to the case of lc-surfaces}
\label{sec:new-adapt-argum}

Let us see how the arguments presented in \S \ref{s:walkthrough} and \S \ref{s:howFR} can be adapted to
the setting of lc-surfaces and used to conclude to the proof of \eqref{e:dream}.  First, we
generalise the upper bound on the number of local types on tangle-free surfaces, Theorem
\ref{cor:TF_curves}, in the following statement.
 
\begin{prp}[{\cite[Theorem A.4]{Moebius-paper}}]
  \label{p:count_with_tangles_bis}
  If $L=A\log g$, $\tfL=\kappa \log g$, if $|\Sf|=0$ and if $\chi(\Sf)< \chictau$,
$$\sharp \LocTFlc = \O[ \chictau, \kappa, A]{ 
  g^{10 \chitau\chictau \kappa} }.$$
\end{prp}
Note that, as opposed to Theorem \ref{cor:TF_curves}, this upper bound is not logarithmic in
the genus. This is due to the presence of factors that are exponential in the tangle-free length
$\omega = \kappa \log(g)$. The important thing here is that, for any $\eta>0$, we can pick $\kappa$
small enough so that the quantity grows slower than $g^\eta$.

Then, in \cite{Ours2}, we show that Theorems \ref{thm:exist_asympt_type_intro} and \ref{t:leviathan}
admit the following generalization.

\begin{nota}
  Let us fix a function $\nu$ on $\mathbf{M}$, assumed to satisfy the following properties.
  \begin{itemize}
  \item If $Z \in\mathcal{M}_{(\vg, \vn)}^*$ has 1d-part $c=(c_1, \ldots, c_j)$ and 2d-part
    $\sigma=(\sigma_1, \ldots, \sigma_m)$, and if $\nu(Z)\not=0$, then
    \begin{itemize}
    \item $\ell(c_i)\leq \kappa$ for $i=1, \ldots, j$;
    \item $\ell^{\mathrm{max}}(\partial \sigma_i)\leq 3\tfL \chi(\sigma_i)$ for $i=1, \ldots, m$.
    \end{itemize}
  \item There exists a bound of the form $|\nu(Z)| \leq \mathfrak{f}(|Z|, \chi(Z))$, i.e. a bound on
    $\abso{\nu(Z)}$ depending only on $|Z|=j$ and on $\chi(Z)=\sum_{i=1}^m \chi(\sigma_i)$. We can
    assume without loss of generality that
    $\chi \mapsto \mathfrak{f}(j,\chi)$ is a non-decreasing function for any fixed $j$.
  \end{itemize}
\end{nota}

\begin{rem} The hypotheses above are satisfied for $\nu= \mu_{\kappa, \tfL}$, the Moebius function
  constructed in Theorem \ref{t:moebius}. Indeed, $\mu_{\kappa, \tfL}$ satisfies the last hypothesis
  with the bound
  \begin{equation}
    \label{eq:bound_nu_mu}
    \mathfrak{f}(j, \chi)  := \frac1{j!} U(\chi) \, e^{\tfL V(\chi)}
  \end{equation}
  where $U$ and $V$ are the sequences appearing in \eqref{e:uppermu}.
\end{rem}

\begin{thm}[{\cite{Ours2}}]
  \label{thm:leviathan_lc}
   Let ${\type}=\eqc{{{\Sf, \curve, \tau}}}$ be a local lc-surface in $\LocTFlc$ such that
   $|\Sf|=0$.
   \begin{enumerate}
   \item For any $j\geq 0$, there exists a unique family of continuous functions
     $( h_k^{\nu, {\type}, j})_{k \geq 0}$ such that, for any test function $F$ and any integer
     $K \geq 0$, any large enough $g$, any $\eta>0$,
  \begin{align*}
    \avN[R_j{\type}]{\nu \otimes F} = \sum_{k=\chi(\Sf)}^K
    &
      \frac{1}{g^k} \int_0^{+ \infty} F(\ell) h_k^{\nu, {\type}, j}(\ell) \d \ell \\
    & +  \O[\Sf,K, \eta, \chitau, Q]{ e^{6 \chitau \tfL }  a(\kappa, K)^j\,  \mathfrak{f}(j,
      \chitau)
      \frac{\norm{ F(\ell) \, e^{(1+\eta)\ell}}_\infty}{g^{K+1}}}
  \end{align*}
  for some explicit $a(\kappa, K)$.
\item For any $k$, $\ell \mapsto h_k^{\nu, {\type}, j}(\ell) $ is a Friedman--Ramanujan
  function in the weak sense.  More precisely, there exists $m(k, \Sf)$ and $c(k, \Sf)$, depending only
  on the filled surface $\Sf$, and not on $(\curve, \tau)$ nor on the integer $j$, such that for all
  $j$, $ h_k^{\nu, {\type}, j}\in \cF_w^{m(k, \Sf), c(k, \Sf)}$, and
  \begin{equation}
    \label{eq:bound_FR_norm_lc}
    \norm{ h_k^{\nu, {\type}, j} }_{\cF_w^{m(k, \Sf), c(k, \Sf)}}
    \leq c(k, \Sf) \, e^{6 \chitau \tfL} a(\kappa, K)^j\,   \mathfrak{f}(j, \chitau).
  \end{equation}
\end{enumerate}
 \end{thm}

 The decay in $1/j! $ in \eqref{eq:bound_nu_mu} is particularly useful to compensate with $ a(\kappa, K)^j$ when summing the previous
 result over all values of $j$. Indeed, we can use it to write, for $|\Sf|=0$ and
 ${\type}\in \LocTFlc$,
    \begin{align*}
      \sum_{j=0}^{+\infty}   \avN[R_j{\type}]{\mu_{\kappa, \tfL} \otimes F}=
      & \sum_{k=\chi(\Sf)}^K 
        \frac{1}{g^k} \int_0^{+ \infty} F(\ell) h_k^{\mu_{\kappa,\tfL},{\type}}(\ell) \d \ell
      \\ &  +   \O[\Sf,K, \eta, \chitau, Q]{e^{6 \chitau \tfL }   \exp(a(\kappa, K) ) \,  U(\chitau) \,
           e^{\tfL V(\chitau)} \frac{\norm{ F(\ell) \, e^{(1+\eta)\ell}}_\infty}{g^{K+1}}}
    \end{align*}
    where $h_k^{\nu, {\type}} := \sum_{j=0}^{+\infty}h_k^{\nu, {\type}, j}$ is a convergent
    sum. Then, Theorem \ref{thm:leviathan_lc} tells us that $h_k^{\nu, {\type}}$ is a weak Friedman--Ramanujan function,
    and that there exists $m(k, \Sf)$ and $c(k, \Sf)$ such that for all $k$,
    $h_k^{\mu_{\kappa,\tfL},{\type}}\in \cF_w^{m(k, \Sf), c(k, \Sf)}$, with
    $$\norm{h_k^{\mu_{\kappa,\tfL},{\type}}}_{\cF^{m(k, \Sf), c(k, \Sf)}} \leq c(k, \Sf) \, e^{6
      \chitau \tfL}   \exp(a(\kappa, K) ) \,  U(\chitau) \, e^{\tfL V(\chitau)}.$$

    We can now use the Friedman--Ramanujan property to perform integration by parts, again relying
    on Proposition \ref {cor:FR_implies_small}.

 \begin{prp} 
  \label{cor:FR_implies_small_final} Let $\Sf$ be a filling type such that $|\Sf|=0$.
  Let ${\type}\in \LocTFlc$. For any integer $K \geq 0$, there exists
  constants $c(K, \Sf), m(K, \Sf) \geq 0$ (depending only on the filled surface $\Sf$) such that for any large enough
  $g$, any $m \geq m(K, \Sf)$, any $\eta >0$ and $L \geq 1$, 
  \begin{multline*}\label{e:contribT}
    \sum_{j=0}^{+\infty}   \BigavN[R_j{\type}]{\mu_{\kappa, \tfL} \otimes (\ell \, e^{- \frac \ell 2} \, \D^m H_L(\ell))} 
   \\ =  \O[m,\Sf,K,\eta]{e^{6 \chitau \tfL }   \exp(a(\kappa, K) ) \,  U(\chitau) e^{\tfL V(\chitau)}
      \left(L^{c_K^\Sf} + \frac{e^{L(\frac 1 2  +\eta)}}{g^{K+1}}\right)}.
  \end{multline*}
\end{prp}
Finally summing over all ${\type}\in \LocTFlc$ in \eqref{e:totalj}, we obtain an integer $m$ such
that for any filling type $\Sf$ such that $|\Sf|=0$ and $\chi(\Sf) \leq \chictau$, 
\begin{equation}
  \label{c:almost}
  \begin{split}
    &\Big| \sum_{j=0}^{+\infty}   \sum_{{\type}\in \LocTFlc}
       \BigavN[R_j{\type}]{\mu_{\kappa, \tfL} \otimes  (\ell \, e^{- \frac \ell 2} \, \D^m H_L(\ell))}\Big| \\
    &  =        \O[m,\Sf,K,\eta]{                                                                                                    \sharp \LocTFlc  e^{6\chitau\tfL }  \exp(a(\kappa, K) ) \,  U(\chitau) e^{\tfL V(\chitau)}
                                                                                                     \left(L^{c_K^\Sf} + \frac{e^{L(\frac 1 2  +\eta)}}{g^{K+1}}\right)},
  \end{split}
\end{equation}
which will be exactly what we need in order to conclude, in light of Proposition
\ref{p:count_with_tangles_bis} and equations \eqref{e:conditioned} and \eqref{e:muav}.

\subsection{Final choice of parameters, and conclusion}

Let $\alpha>0$ and $0<\eps< \frac14 -\alpha^2$.
We pick our parameters the following way.
\begin{itemize}
\item The order $K$ of the expansion is chosen so that $\frac{1}{2(K+2)} < \alpha$, slightly strengthening the previous constraint $\frac{1}{2(K+2)} \leq \alpha$ in order to leave room for additional polynomial losses in $g$ (see the choice of $\kappa$ below).
\item $L=A\log g$ with $A=2(K+2)$; as before, this ensures that $\frac{e^{L/2}}{g^{K+1}}=g$, to
  match the topological term of the trace formula.
\item $\chic = 3A$ in order to have \eqref{e:conditioned}.
\item $\chitau=2A+4$ and $Q=2A+5+6\chitau$, as recommended in Proposition \ref{prp:proba_small_Moebius}.
\item $\chictau=\chictau(A, \chitau, Q)$ is chosen to have \eqref{e:chi_reduc}, so that we can
  restrict attention to c-filling types $\Sf$ with $\chi(\Sf)< \chictau$.
\item The integer $m$ in the power $\cD^m$ is chosen as the maximum of the integers given by
  Proposition \ref{cor:FR_implies_small_final}, for the finite number of c-filling types $\Sf$ with $\chi(\Sf) < \chictau$ and $|\Sf|=0$.
\item $\tfL=\kappa \log g$, where $0<\kappa< \min(2/3, 2\argsh 1)$ is chosen so that $2\kappa   V(\chitau)+12\chitau\kappa <1$ as recommended in Proposition \ref{prp:proba_small_Moebius}, and also
 $10\chitau \chictau\kappa+ 6\chitau\kappa+ 2\kappa   V(\chitau)< \alpha A -1$, to ensure that the upper bound obtained in equation \eqref{c:almost} is $\cO(e^{(\alpha+\eta)L})$.
\item In Proposition \ref{cor:FR_implies_small_final}, we take $\eta<\eps $.
\end{itemize}

Then, for $F(\ell)= \ell \cD^m H_L(\ell)e^{-\ell/2}$, all terms in \eqref{e:muav} are
$o(e^{(\alpha + \eps)L})$. By \eqref{e:conditioned}, this allows to conclude that
 \begin{align*} 
 \limsup   \Pwp{\delta \leq \lambda_1 \leq \frac 1 4 - \alpha^2 - \epsilon} \leq  
 \kappa^2 .
    \end{align*}
Since $\kappa$ could be arbitrarily small, this proves \eqref{e:dream}.

\bibliographystyle{plain}
\bibliography{bibliography}

 \end{document}